\def \NN {\mathbb N}
\def \CC {\mathbb C}

\def \RR {\mathbb R}
\def \ZZ {\mathbb Z}

\def \epsilon{\varepsilon}

\def \A  {{\mathcal A}}

\def \D  {{\mathcal D}}

\def \LL {{\mathcal L}}

\def \R  {{\mathcal R}}

\def \d {\text{d}}

\def \bfalpha {{\boldsymbol{\alpha}}}

\def \fine {{\hfill \qedsymbol}}

\def \si {\sigma}

\newcommand{\res}{\text{res}}
\renewcommand{\S}{{\mathcal S}}

\documentclass[12pt,reqno]{amsart}

\usepackage{amsmath,amssymb}
\numberwithin{equation}{section}

\begin{document}

\title[]{Twists and resonance of $L$-functions, I}

\author[]{J.Kaczorowski  \lowercase{and} A.Perelli}
\date{}
\maketitle

\bigskip
{\bf Abstract.} We obtain the basic analytic properties, i.e. meromorphic continuation, polar structure and bounds for the order of growth, of all the nonlinear twists with exponents $\leq 1/d$ of the $L$-functions of any degree $d\geq 1$ in the extended Selberg class. In particular, this solves the resonance problem in all such cases.



\vskip1cm
\section{Introduction}

\smallskip
{\bf Statement of results.} 
Given an $L$-function $F(s)$ of positive degree $d$ in the extended Selberg class $\S^\sharp$ (namely, the class of Dirichlet series with meromorphic continuation and functional equation; see below for definitions and notation), in \cite{Ka-Pe/2005} we considered the {\it standard twist}
\begin{equation}
\label{1-1}
F(s,\alpha) = \sum_{n=1}^\infty \frac{a(n)e(-\alpha n^{1/d})}{n^s} \hskip1.5cm \alpha\in\RR, \ \alpha\neq0, \ e(x) = e^{2\pi i x}
\end{equation}
and obtained its main analytic properties; see Theorems 1 and 2 in \cite{Ka-Pe/2005}. Precisely, writing 
\[
\text{Spec}(F) = \{\alpha>0: a(n_\alpha)\neq 0\} \ \text{where} \ n_\alpha = qd^{-d}\alpha^d \ \text{and} \ a(n_\alpha) = 
\begin{cases}
0 \ & \text{if} \ n_\alpha\not\in\NN \\
a(n_\alpha) & \text{if} \ n_\alpha\in\NN,
\end{cases}
\]
(here $q=q_F$ is the conductor of $F(s)$, see below) we proved that $F(s,\alpha)$ is an entire function if $|\alpha|\not\in \text{Spec}(F)$, while $F(s,\alpha)$ is meromorphic over $\CC$ if $|\alpha|\in \text{Spec}(F)$. In the latter case $F(s,\alpha)$ has at most simple poles at the points
\begin{equation}
\label{1-2}
s_k =\frac{d+1}{2d} -\frac{k}{d} -i\theta_F \hskip2cm k=0,1,...,
\end{equation}
and
\begin{equation}
\label{1-3}
\res_{s=s_0} F(s,\alpha) = \frac{c_0(F)}{q^{s_0}} \frac{\overline{a(n_{|\alpha|})}}{n_{|\alpha|}^{1-s_0}},
\end{equation}
where $c_0(F)\neq 0$ is a certain constant depending only on $F(s)$, and $\theta_F$ is the internal shift of $F(s)$ (see below, and note a slight change in the definition of $\theta_F$ with respect to our previous papers). Moreover, $F(s,\alpha)$ has polynomial growth on vertical strips.

\medskip
{\bf Remark 1.} 
Actually, in \cite{Ka-Pe/2005} we considered only the case with $\alpha>0$, but $F(s,-\alpha) = \overline{\overline{F}(\overline{s},\alpha)}$ (see below for notation) and hence the case with $\alpha<0$ follows at once, since $\theta_{\overline{F}}=-\theta_F$ and $c_0(\overline{F})=\overline{c_0(F)}$. \fine

\medskip
The standard twist $F(s,\alpha)$ plays a relevant role in the Selberg class theory and, moreover, is a new object in the theory of classical $L$-functions (i.e. the $L$-functions associated with algebraic, geometric and automorphic structures). For example, the properties of the standard twist were crucial in the classification of the degree 1 functions obtained in \cite{Ka-Pe/1999a} and in the proof of the degree conjecture for $1<d<2$ devised in \cite{Ka-Pe/2011a}. In the latter paper we also studied the general {\it nonlinear twist}
\begin{equation}
\label{1-4}
F(s;f) = \sum_{n=1}^\infty \frac{a(n)e(-f(n,\bfalpha))}{n^s}
\end{equation}
with functions $f(n,\bfalpha)$ of type
\begin{equation}
\label{1-5}
f(n,\bfalpha) = \sum_{j=0}^N \alpha_j n^{\kappa_j} \hskip1.5cm 0<\kappa_N<\dots<\kappa_1< \kappa_0, \ \alpha_j\in\RR,
\end{equation}
in the case $\kappa_0>1/d$ and $\alpha_0>0$. In \cite{Ka-Pe/2011a} we could not get a full description of the analytic properties of such twists, but we obtained a useful transformation formula relating $F(s;f)$ to its dual twist $\overline{F}(s;-f^*)$. See Theorem 1.1 in \cite{Ka-Pe/2011a} for the precise statement and below for notation; again, note a slight difference of notation with respect to \cite{Ka-Pe/2011a}, this time concerning $\overline{F}(s;f)$.

\medskip
The properties we prove in this paper for twists of type \eqref{1-4} allow to solve, with a certain degree of generality, the {\it resonance problem} for the $L$-functions from the class $\S^\sharp$. Indeed, the resonance problem may be stated in general terms as describing under what circumstances the coefficients of an $L$-function, once suitably twisted and averaged, show cancellation (no resonance case) or an asymptotic behavior (resonance case). In an essentially equivalent form, the resonance problem may be stated as describing the evolution of meromorphic continuation, polar structure and order of growth of an $L$-function under a suitable set of twists. This is the form of the problem we deal with in this paper, and it is well known that results of the previous form can be deduced by standard methods. We shall discuss the problem in greater detail after the statement of our results.

\medskip
In this paper we study the nonlinear twist \eqref{1-4} for functions of type \eqref{1-5} in the remaining cases, i.e. when $\kappa_0=1/d$ and there is at least one $\alpha_j\neq0$. For completeness, we also study the simpler case where all the exponents $\kappa_j$ are negative. It turns out that a technique different from the one we used in \cite{Ka-Pe/2005} allows a rather complete description of the analytic properties of $F(s;f)$ in these cases. Such a technique is closer to the one we made use of in \cite{Ka-Pe/twist} to deal with the linear twists in degree $d=2$. Moreover, our present results provide further information in the special case \eqref{1-1} as well. We start with the new theorems in such a special case, namely with general bounds on the growth of $F(s,\alpha)$ for $s\in\CC$ (not only for $s$ in vertical strips as in \cite{Ka-Pe/2005}),
since these results will be needed in the proof of the general case. For completeness, in the statements of Theorems 1 and 2 below we also recall some of the properties already established in \cite{Ka-Pe/2005}. Since we deal with functions of degree $d>0$, and it is known that there are no functions in $\S^\sharp$ with degree $0<d<1$ (see Conrey-Ghosh \cite{Co-Gh/1993}, and \cite{Ka-Pe/1999a}), we may assume that $d\geq 1$. Denoting by $\S^\sharp_d$ the set of functions in $\S^\sharp$ with degree $d$, we have the following results.

\medskip
{\bf Theorem 1.} {\sl Let $F\in\S^\sharp_d$ with $d\geq 1$ and let $\alpha\neq0$, $|\alpha|\notin {\rm Spec}(F)$. Then $F(s,\alpha)$ is entire of order $\leq 1$. Moreover, for every $0<\delta<1$ there exist $A,B,C>0$, depending on $F(s)$ and $\delta$, such that for every $s\in\CC$
\[
F(s,\alpha) \ll A^{|\sigma|} (1+|s|)^{d|\sigma|/\delta + B} e^{C|s|^{\delta'}},
\]
where $\delta' = \max(0,2-1/\delta)$ and the $\ll$-symbol may depend on $F(s)$, $\alpha$ and $\delta$.}

\medskip
{\bf Theorem 2.} {\sl Let $F\in\S^\sharp_d$ with $d\geq 1$ and let $\alpha\neq0$, $|\alpha|\in {\rm Spec}(F)$. Then $F(s,\alpha)$ is meromorphic on $\CC$ with at most simple poles at the points $s_k$ in \eqref{1-2}, and residue at $s_0$ given by \eqref{1-3} for $\alpha>0$. Moreover, for every $\delta_0>0$, $0<\delta<1$ and $\eta>1/\delta$ there exist $A,B,C>0$, depending on $F(s)$, $\delta_0$, $\delta$ and $\eta$, such that for $|s-s_k|\geq \delta_0$
\[
F(s,\alpha) \ll A^{|\sigma|} (1+|s|)^{\eta d|\sigma| + B} e^{C(|s|^{\delta'}+|s|^{(3-\eta\delta)/2})},
\]
where $\delta' = \max(0,2-1/\delta)$ and the $\ll$-symbol may depend on $F(s)$, $\alpha$, $\delta_0$, $\delta$ and $\eta$.}

\medskip
{\bf Remark 2.} Note that, choosing e.g. $\delta=1/2$ and $\eta=6$, the bounds in Theorems 1 and 2 imply that $F(s, \alpha)$ has polynomial growth on vertical strips. However, the resulting bounds are essentially qualitative, although explicit values can be obtained for the constant $B$. In this paper we are not looking for sharp results in this respect; for example, one can immediately get sharper bounds by convexity. We refer to Theorem 2 of \cite{Ka-Pe/2005} for $\alpha$-uniform polynomial bounds, again quantitatively not sharp. \fine

\medskip
The proof of Theorems 1 and 2 forms the bulk of the paper. Indeed, such theorems together with Remark 5 below contain the basic bounds on the growth of $F(s,\alpha)$ needed to trigger a kind of iterative process, see Theorem 5 below, leading to the following general result.

\medskip
{\bf Theorem 3.} {\sl Let $F\in\S^\sharp_d$ with $d\geq 1$ and let $f(n,\bfalpha)$ be as in \eqref{1-5} with $\kappa_0=1/d$ and $\alpha_1\neq 0$. Then the nonlinear twist $F(s;f)$ in \eqref{1-4} is entire of order $\leq 1$. Moreover, there exist $A,B,C>0$ and $0\leq \delta<1$, depending on $F(s)$, such that for every $s\in\CC$
\[
F(s;f) \ll A^{|\sigma|} (1+|s|)^{|\sigma|/\kappa_1 + B} e^{C|s|^{\delta}},
\]
where the $\ll$-symbol may depend on $F(s)$ and $f(n,\bfalpha)$.}

\medskip
{\bf Remark 3.} At present we cannot prove, in general, that the twists $F(s;f)$ in Theorem 3 have polynomial growth on vertical lines; see the discussion below about the resonance problem for further information on this point. \fine

\medskip
Theorem 3 is a consequence of Theorems 1, 2 and 5, and its proof is given after the statement of Theorem 5. For completeness we state and prove a simpler and more general result concerning the nonlinear twists with negative exponents. We present this result in a form suitable to be coupled with the previous results by an iterative argument (see Remark 7 below). Let $F(s)$ be an absolutely convergent Dirichlet series for $\sigma>1$, with meromorphic continuation to $\CC$ and singularities contained in a horizontal strip of finite height. Moreover, let 
\[
f(n,\bfalpha)=\sum_{j=0}^N\alpha_jn^{-\kappa_j} \hskip2cm 0<\kappa_0<\dots<\kappa_N, \ \alpha_j\in\RR. 
\]
Denote by $\mu_F(\sigma)$ (resp. $\mu_F(\sigma;f)$), if it exists, the Lindel\"of $\mu$-function of $F(s)$ (resp. $F(s;f)$).

\medskip
{\bf Theorem 4.} {\sl Let $F(s)$ and $f(n,\bfalpha)$ be as above. Then the twist $F(s;f)$ in \eqref{1-4} is meromorphic on $\CC$ with singularities in the same horizontal strip as $F(s)$, and $F(s;f)$ is entire if $F(s)$ is entire. Moreover
\[
\mu_F(\sigma;f) = \mu_F(\sigma)
\]
for $\si$ in any right half-line where $\mu_F(\si)$ exists}.

\medskip
{\bf Remark 4.} From the proof of Theorem 4 it is easy to detect the location of the poles of $F(s;f)$ from the location of the poles of $F(s)$, see \eqref{T4-1} in Section 4. \fine

\medskip
In order to state the results allowing the iterative process leading to Theorem 3, we need to introduce further notation. For $\rho\geq 1$ and $\tau\geq 0$ let $M(\rho,\tau)$ be the class of Dirichlet series $F(s)$, absolutely convergent for $\sigma>1$, admitting holomorphic continuation to $|t|\geq \tau$ and for which there exist $A,B,C>0$ and $0\leq \delta<1$, all depending on $F(s)$, such that for $|t|\geq \tau$
\begin{equation}
\label{1-6}
F(s) \ll A^{|\sigma|} (1+|s|)^{\rho|\sigma|+B} e^{C|s|^\delta},
\end{equation}
where the $\ll$-symbol may depend on $F(s)$, $\rho$ and $\tau$. Moreover, we denote simply by $M(\rho)$ the class $M(\rho,0)$ and write for $\lambda>0$ and $\alpha\in\RR$
\[
F^\lambda(s,\alpha) = \sum_{n=1}^\infty \frac{a(n) e(-\alpha n^\lambda)}{n^s};
\]
thus in particular $F(s,\alpha) = F^{1/d}(s,\alpha)$. The relevance of the class $M(\rho,\tau)$ is clarified by Theorems 1 and 2, and by the following

\medskip
{\bf Remark 5.} If $F\in\S_d^\sharp$ then $F\in M(d, \tau)$ for every $\tau>0$. In fact, $F(s)$ is holomorphic except possibly at $s=1$ (thus the entire functions $F\in\S^\sharp_d$ belong to $M(d)$). Moreover $F(s)$ is bounded for $\sigma\geq 2$, is of polynomial growth (in particular) for $-1 \leq \sigma \leq 2$ and for $\sigma\leq -1$ satisfies
\[
F(s) \ll A^{|\sigma|} (1+|s|)^{d(|\sigma|+1/2)},
\]
thanks to the functional equation and Stirling's formula; see e.g. Lemma 2.1 of \cite{Ka-Pe/twist}. \fine

\medskip
{\bf Theorem 5.} {\sl If $F\in M(\rho,\tau)$ with $\tau>0$, and if $0<\lambda<1/\rho$ and $\alpha\neq 0$, then $F^\lambda(s,\alpha)$ belongs to $M(1/\lambda)$. Moreover, if $F\in M(\rho)$ and $0<\lambda<1/\rho$, then $F^\lambda(s,\alpha)$ belongs to $M(\rho)$ for every $\alpha\in\RR$.}

\medskip
{\bf Remark 6.} When $F\in M(\rho)$ and $0<\lambda<1/\rho$, the entire function $F^\lambda(s,\alpha)$ has the representation
\[
F^\lambda(s,\alpha) = \sum_{k=0}^\infty \frac{(-1)^k}{k!} F(s-\lambda k) (2\pi i\alpha)^k,
\]
the series being absolutely and uniformly convergent for $s$ in compact subsets of $\CC$; see the proof of Theorem 5. \fine

\medskip
{\bf Proof of Theorem 3.} As we already noticed, Theorem 3 is a direct consequence of Theorems 1, 2 and 5. Indeed, if $\alpha_0\neq0$ we start with $F(s,\alpha_0)$. Then we choose $\delta$ (and $\eta$) in Theorems 1 or 2, depending on $|\alpha_0|\not\in$ Spec$(F)$ or $|\alpha_0|\in$ Spec$(F)$, in such a way that $F(s,\alpha_0)$ belongs to $M(\rho,\tau)$ for some $\tau>0$ and $d<\rho<1/\kappa_1$. If $\alpha_0=0$ we use Remark 5, asserting that $F\in M(d,\tau)$ for every $\tau>0$. Thus, for any $\alpha_0\in\RR$ we have that $F(s,\alpha_0)$ belongs to $M(\rho,\tau)$ with some $\tau>0$ and $d\leq\rho<1/\kappa_1$. Then we apply the first part of Theorem 5 with $\lambda=\kappa_1$ to $F(s,\alpha_0)$, thus obtaining that
\[
\sum_{n=1}^\infty \frac{a(n)}{n^s} e(-\alpha_0 n^{1/d} - \alpha_1 n^{\kappa_1})
\]
belongs to $M(1/\kappa_1)$. Now we apply iteratively the second part of Theorem 5 with $\lambda=\kappa_2, \kappa_3,...$, and Theorem 3 follows. \fine

\medskip
{\bf Remark 7.} Thanks to Theorem 4, a similar, but simpler, iterative argument may be applied to any of the twists $F(s,\alpha)$ or $F(s;f)$ considered in Theorems 1, 2 and 3, thus getting meromorphic continuation and polar structure (see Remark 4) of any twist of type \eqref{1-4} with 
\[
f(n,\bfalpha) = \sum_{j=0}^N \alpha_j n^{\kappa_j}
\]
with $\alpha_j\in\RR$ and $\kappa_N<\dots<\kappa_1<\kappa_0=1/d$. \fine

\medskip
{\bf Remark 8.} We finally remark that, as far as we know, the results in Theorems 1, 2 and 3 are new also in the case of classical $L$-functions. It is interesting to note how the behavior in $s$ of the involved functions, essentially as $\sigma\to-\infty$, appears to be critical in order to deduce the properties of their nonlinear twists. This complements the well known importance of the behavior of $L$-functions on vertical strips. A similar phenomenon already arises in \cite{Ka-Pe/twist}, where the behavior of the linear twists as $\sigma\to-\infty$ is shown to give control on the shape of the Euler product; see Theorem 1 of  \cite{Ka-Pe/twist}. \fine

\bigskip
{\bf The resonance problem.}
Theorems 1, 2 and 3 provide the basic analytic properties of the nonlinear twists $F(s;f)$ in all cases where $f(n,\bfalpha)$ has positive exponents $\leq 1/d$. Actually, the same holds if all exponents are $\leq 1/d$ thanks to Theorem 4 (see Remark 7), but in what follows we restrict for simplicity to positive exponents. In particular, polar structure and order of growth of $F(s;f)$ are determined by the above theorems for any function $F\in\S^\sharp_d$ in the following form:

\smallskip
{\sl i)} if $f(n,\bfalpha)= \alpha_0n^{1/d}$ with $|\alpha_0|\in$ Spec$(F)$, then $F(s;f)$ has a simple pole at $s=s_0$, possible simple poles at points $s=s_k$ with $k\geq 1$ (see \eqref{1-2}) and polynomial growth on vertical lines;

{\sl ii)} if $f(n,\bfalpha) = \alpha_0n^{1/d}$ with $0\neq|\alpha_0|\not\in$ Spec$(F)$, then $F(s;f)$ is entire and has polynomial growth on vertical lines;

{\sl iii)} if $f(n,\bfalpha)= \alpha_0n^{1/d} + \alpha_1 n^{\kappa_1}+\dots$ with any $\alpha_0\in\RR$, $\alpha_1\neq 0$ and $\kappa_1<1/d$, then $F(s;f)$ is entire, and on any fixed vertical strip
\begin{equation}
\label{1-7}
F(s;f) \ll \exp(|t|^\delta)
\end{equation}
with some $\delta<1$.

\smallskip
Thus such results solve the resonance problem, in the second form stated above, for all nonlinear twists with exponents $\leq 1/d$ of any function of degree $d\geq 1$ from $\S^\sharp$. Moreover, in cases {\sl i)} and {\sl ii)} standard techniques can be used to describe the behavior of the smoothed nonlinear exponential sums
\begin{equation}
\label{1-8}
S_F(x;f,\phi) = \sum_{n=1}^\infty a(n)e(-f(n,\bfalpha)) \phi(\frac{n}{x}),
\end{equation}
where $\phi(u)$ is a smooth function on $(0,\infty)$ with compact support. Precisely, we have
\begin{equation}
\label{1-9}
S_F(x;f,\phi) = \sum_{k\leq K} c_k(F;f)\widetilde{\phi}(s_k) x^{s_k} +O(x^{\frac{d-1}{2d} - \frac{K}{d}}) \hskip1.5cm x\to\infty,
\end{equation}
for any fixed $K\geq 0$, where $\widetilde{\phi}(s)$ is the Mellin transform of $\phi(\xi)$, $c_0(F;f)\neq 0$ in case {\sl i)} and $c_k(F;f)=0$ for $k\geq 0$ in case {\sl ii)}. Interesting applications would follow from suitable uniform bounds in $\bfalpha$, but at present the quality of such bounds is definitely weak; see e.g.  Theorem 2 in \cite{Ka-Pe/2005}. Formula \eqref{1-9} follows from the estimate
\begin{equation}
\label{1-10}
 \widetilde{\phi}(s) \ll |t|^{-h} \hskip1.5cm |t|\to\infty
 \end{equation}
for every fixed $h>0$, uniformly on any fixed vertical strip. In turn, \eqref{1-10} can be obtained from the definition of $\widetilde{\phi}(s)$ by repeated partial integrations:
 \[
 \widetilde{\phi}(s) = \frac{(-1)^h}{s(s+1)\cdots(s+h-1)} \int_0^\infty\phi^{(h)}(u) u^{s+h-1} \d u.
 \]
 In case {\sl iii)}, when only \eqref{1-7} is available, \eqref{1-10} is too weak to deal with the sums \eqref{1-8}. However, in this case we can describe, again by standard techniques, the behavior of smoothed nonlinear exponential sums of type
\[
S_F(x;f,r) = \sum_{n=1}^\infty a(n)e(-f(n,\bfalpha)) e^{-(n/x)^r}
\]
with $r>0$ arbitrary. In this case, for any fixed $K\geq 0$ we have
\[
S_F(x;f,r) = \sum_{k\leq K}\frac{(-1)^k}{k!}F(-kr;f)x^{-kr} + O(x^{-(K+1)r})  \hskip1.5cm x\to\infty.
\]
This is due to the fact that the Mellin transform of $e^{-u^r}$ is $\frac{1}{r}\Gamma(\frac{s}{r})$, and by Stirling's formula 
\[
\Gamma(\frac{s}{r}) \ll e^{-\pi|t|/2r} |t|^{c(a,b)}
\]
uniformly for $|t|\geq 1$ and $a<\sigma<b$.

\medskip
We omit explicit examples for Theorems 1 to 5 since one may easily construct such examples starting with any classical $L$-function. In a forthcoming paper we shall study the resonance properties of the nonlinear twists of type \eqref{1-4} and \eqref{1-5} with $\kappa_0>1/d$.

\bigskip
{\bf Definitions and notation.} 
Given a function $f(s)$ we write $\overline{f}(s) = \overline{f(\overline{s})}$; in paricular, if $f(s)$ is a Dirichlet series then $\overline{f}(s)$ is the Dirichlet series with conjugate coefficients, called the conjugate of $f(s)$. However, when dealing with the twists of a function $F(s)$ we write
\[
\overline{F}(s;f) = \sum_{n=1}^\infty \frac{\overline{a(n)}e(-f(n,\bfalpha))}{n^s},
\]
i.e. $\overline{F}(s;f)$ is the twist of the conjugate of $F(s)$. A completely analogous notation is used in the case of $F^\lambda(s,\alpha)$. A function $F(s)$ belongs to the Selberg class $\S$ if

\smallskip
\noindent
{\sl i)} $F(s)$ is an absolutely convergent Dirichlet series for $\sigma>1$; 

\noindent
{\sl ii)} $(s-1)^mF(s)$ is
an entire function of finite order for some integer $m\geq 0$; 

\noindent
{\sl iii)} $F(s)$ satisfies a functional equation of type $\Phi(s) = \omega\overline{\Phi}(1-s)$,
where $|\omega|=1$ and
\[
\Phi(s) = Q^s\prod_{j=1}^r\Gamma(\lambda_js+\mu_j)F(s)
\]
with $r\geq 0$, $Q>0$, $\lambda_j>0$, $\Re\mu_j\geq 0$;

\noindent
{\sl iv)} the Dirichlet coefficients $a(n)$ of $F(s)$ satisfy $a(n) \ll n^\epsilon$ for every
$\epsilon>0$; 

\noindent
{\sl v)} $\log F(s)$ is a Dirichlet series with coefficients $b(n)$ satisfying $b(n)=0$ unless
$n=p^m$, $m\geq 1$, and $b(n)\ll n^\vartheta$ for some $\vartheta<1/2$. 

\smallskip
\noindent
The extended Selberg class $\S^\sharp$ consists of the non-zero functions satisfying only axioms {\sl i)}, {\sl ii)} and {\sl iii)}. Degree, conductor and $\xi$-invariant of $F\in\S^\sharp$ are defined respectively by
\[
d_F = 2\sum_{j=1}^r\lambda_j, \hskip1cm q_F= (2\pi)^{d_F}Q^2\prod_{j=1}^r\lambda_j^{2\lambda_j}, \hskip1cm \xi_F = 2\sum_{j=1}^r (\mu_j-1/2) = \eta_F + i\theta_Fd_F;
\]
note the slight change of notation for $\theta_F$, the internal shift of $F(s)$, with respect to our previous papers. We refer to Selberg \cite{Sel/1989} and Conrey-Ghosh \cite{Co-Gh/1993}, to our survey papers \cite{Ka-Pe/1999b}, \cite{Kac/2006}, \cite{Per/2005}, \cite{Per/2004}, \cite{Per/2010} and to our forthcoming book \cite{book} for the basic information and results on the Selberg class.

\medskip
In the next sections, the constants $A,B,C>0$ are sufficiently large and may depend on $F(s)$ and on the parameters specified in the theorems. Moreover, such constants may also depend on other parameters appearing in the lemmas below; if this happens, it will be explicitly stated. The same applies to the constants implicit in the $O$- and $\ll$-notation, as well as to the constants $c,c_1,...>0$. In all cases, the value of such constants may not be the same at each occurrence. If $F\in\S^\sharp_d$, $d\geq 1$, and $\alpha>0$ are fixed we write
\[
F_X(s,\alpha) = \sum_{n=1}^\infty \frac{a(n)}{n^s} e^{-z_Xn^{1/d}}
\]
where $X\geq X_0>0$ is an integer, $X_0$ is sufficiently large and may depend on $F(s)$ and other parameters as for the above constants $A,B,C$, and
\[
z_X=2\pi\alpha \omega_X \hskip2cm \omega_X=\frac{1}{X}+i. 
\]
Clearly, $F_X(s,\alpha)$ converges for every $s\in\CC$ and hence $F_X(s,\alpha)$ is an entire function, 
\[
\lim_{X\to\infty} F_X(s,\alpha) = F(s,\alpha) \hskip1.5cm \sigma>1,
\]
and by Mellin's transform we have
\begin{equation}
\label{1-11}
F_X(s,\alpha) = \frac{1}{2\pi i}\int_{(c)} F(s+\frac{w}{d}) \Gamma(w) z_X^{-w}\d w \hskip1.5cm c>\max(0,d(1-\sigma)),
\end{equation}
where $w=u+iv$, $z_X^{-w} = e^{-w\log z_X}$ and $\log z_X$ is meant as principal value. Let $\sigma\leq 3/2$, $c=d(3-\sigma)$,
\[
M_F=\max_{1\leq j\leq r} \frac{(1+|\mu_j|)^2}{\lambda_j}, \hskip2cm K=[3d(2+M_F)] +\frac12
\]
and $V\geq V_0>0$ be a parameter to be chosen in the proofs, where $V_0$ is sufficiently large and may depend on $F(s)$ and other parameters as for the above constants $A,B,C$. Finally, let 
\[
\begin{split}
\LL_{-\infty} &= (d(2-\sigma)-i\infty, d(2-\sigma) +iV] \\
\LL_V&=[d(2-\sigma) +iV,-K+iV] \\
\LL_\infty &= [-K+iV, -K+i\infty).
\end{split}
\]

\bigskip
{\bf Outline of the proofs.}
To prove Theorems 1 and 2 we start with \eqref{1-11}, with the aim of letting $X\to\infty$. The integral over the negative part of the line $\si=c$ has good convergence properties. Hence we deform the integration as in \eqref{new2}, thus avoiding the poles of the integrand and preparing for the use of the functional equation. The first two terms in \eqref{new2} are dealt with by Lemma A in Section 2, which holds for every $\alpha>0$ and whose proof is based on a direct application of Stirling's formula. In the integral over $\LL_\infty$ we apply the functional equation and expand $\overline{F}(1-s-\frac{w}{d})$, thus getting \eqref{2-3}. The integral $I_X(s,y)$ in \eqref{2-3} depends on the data of $F(s)$ and is close to an incomplete hypergeometric function. We refer to Section 2 of \cite{Ka-Pe/2005} for an analysis of such functions, see Theorem 2.1 there. In particular, it turns out that the limit as $X\to\infty$ of such hypergeometric functions is meromorphic in $s$ for every $y>0$; moreover, it has a simple pole at $s=s_0$ (and possibly at the points $s=s_k$ with $k\geq 1$) for $y=d/\beta^{1/d}$ (see \eqref{2-2}), and is holomorphic otherwise. Such a behavior gives rise to the notion of Spec$(F)$, and this clarifies why the treatment of the function $F^{(3)}_X(s)$ in \eqref{new2} is different depending on $|\alpha|\in$ Spec$(F)$ or not. As we already remarked, in this paper we follow a different approach to the study of $I_X(s,y)$, leading to the new estimates in Theorems 1 and 2 required by Theorems 3 and 5. In view of Remark 1 we deal only with the case $\alpha>0$.

\smallskip
When $0<\alpha\not\in$ Spec$(F)$ we have $a(n_\alpha)=0$, and we deal separately with the integrals $I_X(s,y_n)$ in \eqref{2-3} with $n<n_\alpha$ and $n>n_\alpha$. In both cases the treatment is based on a change of the path of integration and on a careful application of the uniform version of Stirling's formula proved in \cite{Ka-Pe/2011b}, see Lemma D in Section 2. This leads to Lemma B in Section 2. In turn, Lemmas A and B allow the use of Vitali's convergence theorem when $X\to\infty$, and Theorem 1 follows; see the proof of Theorem 1 in Section 2.

\smallskip
When $\alpha\in$ Spec$(F)$ we have the additional term $n=n_\alpha$, and the corresponding integral $I_X(s,d/\beta^{1/d})$ requires a deeper analysis. The starting point is Lemma 3.1, where the $\Gamma$-factors coming from the functional equation are glued to a single term plus smaller order terms. Then in a series of lemmas (from 3.2 to 3.8) we study and transform the expression of $I_X(s,d/\beta^{1/d})$ coming from Lemma 3.1, finally getting Lemma E in Section 3, where the remainder term is of the required form. Next we borrow some arguments from our previous treatment of the standard twist in \cite{Ka-Pe/1999a} and \cite{Ka-Pe/2005}, giving the explicit expression in \eqref{3-29} and describing the polar structure of the limit as $X\to\infty$ of the integrals $I_{X,\nu}(s)$ in Lemma E. Theorem 2 follows then from Vitali's convergence theorem as before, and from bounds of the required form for the expression in \eqref{3-29}, away from its poles.

\medskip
We already proved Theorem 3, and the proof of Theorem 4 follows by the Taylor expansion of $e(-f(n,\bfalpha))$ and the good convergence properties of the resulting series. The proof of Theorem 5 is simpler than the proof of Theorems 1 and 2 thanks to the good convergence properties of the integral in \eqref{1-11}, due to the choice $\lambda<1/\rho$. Indeed, when $F(s)$ is entire and belongs to $M(\rho)$ we shift the line of integration in \eqref{1-11} to the left, thus getting a sum over the residues of $\Gamma(w)$, see \eqref{T5-13}. Thanks to $\lambda<1/\rho$, as $X\to\infty$ we obtain an expression of $F^\lambda(s,\alpha)$ as a series, see \eqref{new1}, which is nicely convergent over $\CC$. Moreover, the required bound for $F^\lambda(s,\alpha)$ follows by plugging in \eqref{new1} the bound for $F(s-\lambda k)$. The case of meromorphic $F(s)$ is technically more delicate, but the good convergence properties of the involved integrals, again due to the choice $\lambda<1/\rho$, are decisive in this case as well.

\bigskip
{\bf Acknowledgements.} 
This research was partially supported by the Istituto Nazionale di Alta Matematica and by grant N N201 605940 of the National Science Centre.

\bigskip
\section{Proof of Theorem 1}

\smallskip
We follow the notation at the end of Section 1 and suppose that $F(s)$ and $\alpha$ are as in Theorem 1. In view of Remark 1 of Section 1 we may assume that $\alpha>0$. Thanks to the polynomial growth of $F(s+\frac{w}{d})$ and the decay of $\Gamma(w)$ on vertical strips, and to the  location of their poles, we have
\begin{equation}
\label{new2}
F_X(s,\alpha) =  \frac{1}{2\pi i}\int_{\LL_{-\infty}\cup\LL_V\cup\LL_\infty} F(s+\frac{w}{d}) \Gamma(w) z_X^{-w}\d w = F_X^{(1)}(s) + F_X^{(2)}(s) + F_X^{(3)}(s),
\end{equation}
say. Recalling that $\delta'=\max(0,2-1/\delta)$ we have

\medskip
{\bf Lemma A.} {\sl Let $\sigma\leq 3/2$ and $V\geq V_0(1+|s|)$. Under the hypotheses of Theorem $1$ and with the notation in Section $1$ we have
\[
|F_X^{(1)}(s)| + |F_X^{(2)}(s)| \ll A^{|\sigma|} V^{d|\sigma| + 5d}
\]
uniformly for $X\geq X_0$.}

\medskip
{\bf Lemma B.} {\sl Let $\sigma\leq 3/2$ and $V\geq V_0(1+|s|)^{1/\delta}$. Under the hypotheses of Theorem $1$ and with the notation in Section $1$ we have
\[
|F_X^{(3)}(s)| \ll A^{|\sigma|} V^{d|\sigma| + B} e^{C(|s|^{\delta'}+|s|V^{\delta-1})} 
\]
uniformly for $X\geq X_0$.}

\medskip
Before proving Lemmas A and B we recall Vitali's convergence theorem (see Section 5.21 of Titchmarsh \cite{Tit/1939}), and show that Theorem 1 in an immediate consequence of Lemmas A, B and Vitali's theorem.

\medskip
{\bf Lemma C.} (Vitali's convergence theorem) {\sl Let $f_X(z)$ be a sequence of holomorphic functions on a region $\D$ and let $|f_X(z)| \leq M$ for every $X\geq X_0$ and $z\in\D$. Suppose that $f_X(z)$ tends to a limit, as $X\to\infty$, at a set of points having an accumulation point in $\D$. Then $f_X(z)$ tends to a limit uniformly in any domain $\D'$ whose closure is contained in $\D$, and hence such a limit is holomorphic and bounded by $M$ on $\D'$.}

\medskip
{\bf Proof of Theorem 1.} The result is obvious for $\sigma \geq 5/4$, while for $\sigma< 5/4$ we recall that $F_X(s,\alpha)$ is holomorphic and we apply Lemmas A and B with the choice $V=V_0(1+|s|)^{1/\delta}$ to get
\begin{equation}
\label{2-1}
F_X(s,\alpha) \ll A^{|\sigma|} (1+|s|)^{d|\sigma|/\delta + B} e^{C|s|^{\delta'}}
\end{equation}
uniformly for $X\geq X_0$. Given any $s=\si+it$ with $\si<5/4$ we consider the region 
\[
\D_s=\{z=x+iy\in\CC: \si-1/10<x<3/2, |y|<|t|+1/10\}
\]
and apply Lemma C with the choice $f_X(z)=F_X(z,\alpha)$ and $\D=\D_s$. Thanks to \eqref{2-1} and since $|z|\leq |s|+c$, $|x|\leq |\si|+c$, we may choose 
\[
M = KA^{|\sigma|} (1+|s|)^{d|\sigma|/\delta + B} e^{C|s|^{\delta'}}
\]
with a suitable constant $K>0$. Recalling that $\lim_{X\to\infty} F_X(z,\alpha) = F(z,\alpha)$ for every $z$ with $x>1$, since $s\in\D_s$ we obtain from Lemma C that $F(s,\alpha)$ has holomorphic continuation to the half-plane $\si<5/4$ and
\[
F(s,\alpha) \ll A^{|\sigma|} (1+|s|)^{d|\sigma|/\delta + B} e^{C|s|^{\delta'}}.
\]
Theorem 1 is therefore proved. \fine

\bigskip
{\bf Proof of Lemma A.} We start with the estimation of $F_X^{(1)}(s)$. For $w=u+iv\in\LL_{-\infty}$ we have $u=d(2-\sigma)$, hence $\Re(s+w/d)=2$ and therefore $F(s+w/d)\ll1$. Moreover, by Stirling's formula we have
\[
\Gamma(w) \ll A^{|\sigma|} (1+|w|)^{d(2-\sigma)-1/2} e^{-v\arg w},
\]
and since $v\arg w = v \arctan \frac{v}{d(2-\sigma)} = |v|(\frac{\pi}{2} + O(\frac{1+|\sigma|}{|v|}))$ we get
\[
\Gamma(w) \ll A^{|\sigma|} (1+|w|)^{d(2-\sigma)-1/2} e^{-\pi |v|/2}.
\]
Further we have
\[
|z_X^{-w}| = e^{-d(2-\sigma)\log|z_X| +v\arg z_X} \ll A^{|\sigma|} e^{v(\pi/2 -O(1/X))}.
\] 
Therefore for $V\geq V_0(1+|s|)$ and $X\geq X_0$ we obtain
\[
\begin{split}
F_X^{(1)}(s) &\ll A^{|\sigma|} \int_{-\infty}^V e^{-\pi|v|/2 +v(\pi/2 -O(1/X))} (1+d(2+|\si|)+ |v|)^{d(2+|\sigma|)} \d v \\
&\ll A^{|\sigma|} \left(\int_{-V}^V V^{d(2+|\sigma|)} \d v + \int_{-\infty}^{-V} e^{-|v|(\pi -O(1/X))} (1+|v|)^{d(2+|\sigma|)} \d v
\right) \\
&\ll A^{|\sigma|} V^{d|\si| +2d+1} + A^{|\sigma|}  \Gamma(d(2+|\si|)+1) \ll A^{|\sigma|} (V^{d|\si| +2d+1} + (1+|\si|)^{d|\si|})
\end{split}
\]
uniformly in $X$, and the first assertion of Lemma A follows.

In order to estimate $F_X^{(2)}(s)$ we recall the estimate for $F(s)$ in Remark 5 of the Introduction, which we report in the form
\[
F(s) \ll
\begin{cases}
A^{|\si|} (1+|s|)^{d(|\si|+1/2)} & \text{if} \ \si\leq -1 \\
(1+|s|)^{3d/2} & \text{if} \ \si\geq -1.
\end{cases}
\] 
Accordingly, we split the path of integration into 
\[
[-K+iV, -d(1+\si)+iV] \cup [-d(1+\si)+iV, d(2-\si)+iV],
\]
thus getting
\[
\begin{split}
F_X^{(2)}(s) \ll \ &A^{|\si|} \int_{-K}^{-d(1+\si)} V^{d(|\si+u/d|+1/2)} |\Gamma(u+iV) z_X^{-u-iV}| \d u \\
&+ V^{3d/2} \int_{-d(1+\si)}^{d(2-\si)}  |\Gamma(u+iV) z_X^{-u-iV}| \d u.
\end{split}
\]
Arguing as before and taking into account the definition of $V$ and $K$ we have
\[
|\Gamma(u+iV) z_X^{-u-iV}| \ll A^{|\si|} V^{u-1/2}
\]
uniformly in $X$, hence
\[
F_X^{(2)}(s) \ll A^{|\si|}\left(\int_{-K}^{-d(1+\si)} \hskip-.2cm V^{d(|\si+u/d|+1/2)+u-1/2} \d u + \int_{-d(1+\si)}^{d(2-\si)} \hskip-.2cm V^{(3d + 2u-1)/2} \d u \right) \ll A^{|\si|} V^{d|\si|+5d}
\]
since $\si+u/d<0$ in the first integral, and Lemma A follows. \fine

\medskip
{\bf Remark.} Note that the proof of Lemma A does not depend on the hypothesis that $\alpha\notin$ Spec$(F)$, hence Lemma A holds for any $\alpha>0$. The proof of Lemma B is definitely more delicate, and condition $\alpha\notin$ Spec$(F)$ becomes crucial in such a proof. \fine

\medskip
The following uniform version of Stirling's formula (see \cite{Ka-Pe/2011b}) will be repeatedly used in this paper, hence for convenience we state it as

\medskip
{\bf Lemma D.} {\sl Let $N\geq 0$ be an integer, $D\geq 1$ and let $z,a\in\CC$ satisfy
\[
\Re(z+a)\geq0, \hskip1cm  |a|\leq \frac{3}{5} |z|, \hskip1cm N\leq 2D|z|.
\]
Then}
\[
\begin{split}
\log\Gamma(z+a) &= (z+a-\frac12)\log z - z + \frac12\log2\pi \\
&+ \sum_{j=1}^N \frac{(-1)^{j+1} B_{j+1}(a)}{j(j+1)} \frac{1}{z^j} +O\big(\frac{1}{|z|^{N+1}} \big((N+\frac{|a|^2}{(N+1)^2})|a|^N + D^NN!\big)\big).
\end{split}
\]

\medskip
{\bf Proof of Lemma B.} For $w\in\LL_\infty$ we have $\Re(s+w/d) = \si-K/d \leq 3/2 - 6 +1/2d \leq -4$. Hence, writing
\[
h(w,s) = \prod_{j=1}^r \frac{\Gamma(\lambda_j(1-s) +\overline{\mu}_j - \frac{\lambda_jw}{d})}{\Gamma(\lambda_js +\mu_j + \frac{\lambda_jw}{d})} \Gamma(w),
\]
applying the functional equation and using the Dirichlet series expansion of $\bar{F}(1-s-\frac{w}{d})$ we get
\[
\begin{split}
F^{(3)}_X(s) &= \frac{\omega Q^{1-2s}}{2\pi i} \int_{\LL_\infty} \bar{F}(1-s-\frac{w}{d}) h(w,s)Q^{-2w/d} z_X^{-w} \d w \\
&\ll A^{|\si|} \sum_{n=1}^\infty \frac{|a(n)|}{n^{1-\si}}\left| \int_{\LL_\infty} h(w,s) \big(\frac{Q^{2/d}z_X}{n^{1/d}}\big)^{-w}\d w\right|.
\end{split}
\]
Recalling that $z_X=2\pi\alpha\omega_X$ and $\omega_X=\frac{1}{X}+i$, for convenience we write
\[
I_X(s,y) = \frac{1}{2\pi i}\int_{\LL_\infty} h(w,s) (y\omega_X)^{-w}\d w \hskip1.5cm y>0,
\]
hence recalling the definition of $n_\alpha$ in the Introduction and letting
\begin{equation}
\label{2-2}
\beta=\prod_{j=1}^r \lambda_j^{2\lambda_j} \hskip1.5cm y_n = \big(\frac{n_\alpha}{n}\big)^{1/d} \frac{d}{\beta^{1/d}}
\end{equation}
we have
\begin{equation}
\label{2-3}
F^{(3)}_X(s) \ll A^{|\si|} \sum_{n=1}^\infty \frac{|a(n)|}{n^{1-\si}} |I_X(s,y_n)|.
\end{equation}
Since $\alpha\not\in$ Spec$(F)$ we may assume that $n\neq n_\alpha$ always. In the rest of the proof we obtain suitable bounds for $I_X(s,y_n)$ in the two cases $n>n_\alpha$ and $n<n_\alpha$ which, once inserted in \eqref{2-3}, will prove Lemma B. Actually, in what follows we deal with $I_X(s,y)$ in the general cases $0<y<d/\beta^{1/d}$ and $y>d/\beta^{1/d}$.

\bigskip
{\bf Case $0<y<d/\beta^{1/d}$.} In this case we change the path of integration from $\LL_\infty$ to the half-line $\LL'_\infty$ where $w=u+iV$, $-\infty < u \leq -K$, and consider the 90-degrees sector $\R$ formed by $\LL_\infty$ and $\LL'_\infty$. Thanks to the choice of $K$, the function $h(w,s)$ is holomorphic on $\R$, as well as the function $(y\omega_X)^{-w}$. Moreover, by Cauchy's theorem applied to such sector and the arc $\gamma_R$ of the circle $|w|=R$ inside the sector, as $R\to\infty$ we obtain
\[
I_X(s,y) = \frac{1}{2\pi i}\int_{\LL'_\infty} h(w,s) (y\omega_X)^{-w}\d w.
\]
Indeed, a standard application of the Stirling formula with $|w|=R$ shows that the contribution of the integral over the arc $\gamma_R$ tends to 0 as $R\to\infty$.

\medskip
We first deduce from Lemma D the following useful uniform bound for $|\Gamma(z+a)|$: let $0<\delta<1$ and
\begin{equation}
\label{2-4}
\Re(z+a)\geq 0, \hskip1cm  |a|\leq \min(|z|^\delta, \frac{3}{5}|z|), \hskip1cm |z|\geq 1/2;
\end{equation}
then for $\delta' = \max(0,2-1/\delta)$ we have
\begin{equation}
\label{2-5}
|\Gamma(z+a)| = e^{-\Im(z+a)\arg z-\Re z} |z|^{\Re(z+a)-1/2} e^{O(1+|a|^{\delta'})}.
\end{equation}
Indeed, from Lemma D with $N=0$ and $D=1$ we get
\[
\log|\Gamma(z+a)| = \Re\{(z+a-\frac12)(\log|z| + i\arg z) - z\} + O(\frac{|a|^2}{|z|} + 1)
\]
and \eqref{2-5} follows since the last $O$-term is $O(1+|a|^{\delta'})$.

\medskip
Now we proceed with the estimation of $h(w,s)$ for $w\in\LL'_\infty$. First we use the reflection formula for the $\Gamma$ function to write
\begin{equation}
\label{2-6}
h(w,s) = \pi^{-r} S(w,s) \prod_{j=1}^r \Gamma(\lambda_j(1-s) + \overline{\mu}_j -\frac{\lambda_jw}{d}) \Gamma(1-\lambda_js - \mu_j -\frac{\lambda_jw}{d}) \Gamma(w)
\end{equation}
with
\[
S(w,s) = \prod_{j=1}^r \sin \pi(\lambda_js + \mu_j +\frac{\lambda_jw}{d}).
\]
In order to estimate the first $\Gamma$-factor in the product in \eqref{2-6}, we choose $z=-\frac{\lambda_jw}{d} = \frac{\lambda_j}{d}(|u|-iV)$ and $a=\lambda_j(1-s) +\overline{\mu}_j $ and note that conditions \eqref{2-4} are satisfied thanks to the choice of $K$, provided
\begin{equation}
\label{2-7}
V\geq c(1+\frac{1}{\gamma})^{1/\delta} (1+|s|)^{1/\delta}, \ \ \text{where} \ \gamma = |\log(y\beta^{1/d}/d)|
\end{equation}
and $c=c(F,\delta)>0$ is sufficiently large. Indeed, simple computations show that $|z|\geq \lambda_jK/d \geq 3$
\[
|z| \geq \frac{\lambda_jV}{d} \geq \frac{\lambda_jc}{d}(1+|s|)^{1/\delta} \geq \max(|a|^{1/\delta},\frac{5}{3}|a|),
\]
\[
\Re(z+a) = \frac{\lambda_j|u|}{d} +\lambda_j(1-\si) +\Re\mu_j \geq  \frac{\lambda_jK}{d} - \frac{\lambda_j}{2} -|\mu_j| \geq 3.
\]
Hence we may apply \eqref{2-5}. With the above choice of $z$ and $a$ we have
\[
\Im(z+a) = -\lambda_j(\frac{V}{d} +t) + O(1), \hskip1.5cm \arg z = -\arctan \frac{V}{|u|},
\]
\[
\Re z = \frac{\lambda_j|u|}{d}, \hskip1.5cm \Re(z+a) = \lambda_j(\frac{|u|}{d} + |\sigma|) + O(1),
\]
hence from \eqref{2-5} we obtain
\begin{equation}
\label{2-8}
\Gamma(\lambda_j(1-s) + \overline{\mu}_j -\frac{\lambda_jw}{d}) \ll e^{-\lambda_j(\frac{V}{d}+t)\arctan\frac{V}{|u|} -\lambda_j\frac{|u|}{d}} \big|\frac{\lambda_jw}{d}\big|^{\lambda_j(\frac{|u|}{d} + |\sigma|)+B} e^{C(1+|s|^{\delta'})}.
\end{equation}
To estimate the second $\Gamma$-factor in the product in \eqref{2-6} we choose again $z=-\frac{\lambda_jw}{d}$, and $a=1-\mu_j-\lambda_js$. As before, simple computations show that \eqref{2-5} is applicable and the same formulae hold for the involved quantities. Therefore,
\begin{equation}
\label{2-9}
\Gamma(1-\lambda_js -\mu_j -\frac{\lambda_jw}{d}) \ll e^{-\lambda_j(\frac{V}{d}+t)\arctan\frac{V}{|u|} -\lambda_j\frac{|u|}{d}} \big|\frac{\lambda_jw}{d}\big|^{\lambda_j(\frac{|u|}{d} + |\sigma|)+B} e^{C(1+|s|^{\delta'})}.
\end{equation}
The third $\Gamma$-factor in \eqref{2-6} is estimated by first applying the reflection formula and then \eqref{2-5} (although the standard Stirling formula would suffice here) with the choice $z=|u|-iV$ and $a=1$, thus obtaining
\begin{equation}
\label{2-10}
\Gamma(w) \ll \frac{e^{-\pi V}}{|\Gamma(1+|u|-iV)|} \ll e^{-\pi V +V\arctan\frac{V}{|u|} + |u|} |w|^{-|u|-1/2}.
\end{equation}
Further we have 
\begin{equation}
\label{2-11}
S(w,s) \ll e^{\frac{\pi}{2}d(t+\frac{V}{d})},
\end{equation}
and hence from \eqref{2-6} and \eqref{2-8}-\eqref{2-11}, under \eqref{2-7} and recalling \eqref{2-2}, we obtain
\begin{equation}
\label{2-12}
\begin{split}
h(w,s) &\ll A^{|\si|} e^{-d(\frac{V}{d}+t)\arctan\frac{V}{|u|} -|u|} |w|^{d(\frac{|u|}{d}+|\si|) +B} e^{\frac{\pi}{2}d(t+\frac{V}{d})} \beta^{\frac{|u|}{d}} d^{-|u|} \\
&\hskip1.5cm \times e^{-\pi V +V\arctan\frac{V}{|u|} + |u|} |w|^{-|u|-1/2} e^{C(1+|s|^{\delta'})} \\
&\ll A^{|\si|} e^{dt(\frac{\pi}{2}-\arctan\frac{V}{|u|})} e^{-\frac{\pi}{2}V} |w|^{d|\si|+B} \beta^{\frac{|u|}{d}} d^{-|u|} e^{C(1+|s|^{\delta'})}.
\end{split}
\end{equation}

\medskip
In order to estimate $I_X(s,y)$ we observe that, for $w\in\LL_\infty'$,
\[
|(y\omega_X)^{-w}| \ll |y\omega_X|^{|u|} e^{\frac{\pi}{2}V}
\]
and hence from \eqref{2-12} we get
\begin{equation}
\label{2-13}
\begin{split}
I_X(s,y) &\ll A^{|\si|} e^{C(1+|s|^{\delta'})} \int_K^\infty \big(\frac{\beta^{1/d}}{d}y|\omega_X|\big)^{u} e^{dt(\frac{\pi}{2}-\arctan\frac{V}{u})} |w|^{d|\si|+B} \d u \\
&= A^{|\si|} e^{C(1+|s|^{\delta'})} \left(\int_K^{V^\delta} + \int_{V^\delta}^\infty\right) \big(\frac{\beta^{1/d}}{d}y|\omega_X|\big)^{u} e^{dt(\frac{\pi}{2}-\arctan\frac{V}{u})} |w|^{d|\si|+B} \d u \\
&= I^{(1)}_X(s,y) + I^{(2)}_X(s,y),
\end{split}
\end{equation}
say; note that $\frac{\beta^{1/d}}{d}y<1$ in this case. We deal first with $I^{(1)}_X(s,y)$, observing that for $0<u<V^\delta$ we have $\arctan{V/u} = \pi/2 + O(V^{\delta-1})$. Therefore, if $X>c/\gamma$ with a suitable $c>0$ then $\frac{\beta^{1/d}}{d}y|\omega_X|<1$ and hence
\begin{equation}
\label{2-14}
\begin{split}
I^{(1)}_X(s,y) &\ll A^{|\si|} e^{C(1+|s|^{\delta'}+|s|V^{\delta-1})} V^{d|\si|+B} \int_K^\infty e^{-\gamma_Xu}\d u \\
&\ll A^{|\si|} e^{C(1+|s|^{\delta'}+|s|V^{\delta-1})} V^{d|\si|+B} \frac{y^K}{\gamma_X},
\end{split}
\end{equation}
where $\gamma_X = |\log(\frac{\beta^{1/d}}{d}y|\omega_X|)|$. 

\medskip
Concerning $I^{(2)}_X(s,y)$, again for $X>c/\gamma$  we have
\begin{equation}
\label{2-15}
I^{(2)}_X(s,y) \ll A^{|\si|} e^{C(1+|s|^{\delta'})} e^{d\frac{\pi}{2}|t|} \left(\int_{V^\delta}^\infty e^{-\gamma_X u} u^{d|\si|+B} \d u + V^{d|\si|+B} \int_{V^\delta}^\infty e^{-\gamma_X u} \d u \right),
\end{equation}
and we use the following general result: for $1\leq \xi \leq Y$ we have
\begin{equation}
\label{2-16}
\int_Y^\infty e^{-u}u^\xi \d u \leq e^{-Y}Y^\xi \big(1+\frac{\xi^2}{Y}\big).
\end{equation}
Indeed, integrating by parts twice we have
\[
\int_Y^\infty e^{-u}u^\xi \d u \leq e^{-Y}Y^\xi + \xi e^{-Y}Y^{\xi-1} +\xi(\xi-1) \int_Y^\infty e^{-u}u^\xi \frac{\d u}{u^2}.
\]
But the function $u\to e^{-u}u^\xi$ is decreasing for $u\geq \xi$ and hence
\[
\int_Y^\infty e^{-u}u^\xi \frac{\d u}{u^2} \leq e^{-Y}Y^\xi \int_Y^\infty \frac{\d u}{u^2} = e^{-Y}Y^{\xi-1},
\]
therefore
\[
\int_Y^\infty e^{-u}u^\xi \d u \leq e^{-Y}Y^\xi \big(1 + \frac{\xi}{Y} + \frac{\xi(\xi-1)}{Y}\big)
\]
and \eqref{2-16} follows. We use \eqref{2-16} in the first integral in \eqref{2-15} after the change of variable $\gamma_Xu\to u$, hence with the choice $Y=V^\delta\gamma_X$ and $\xi=d|\si|+B$. Recalling \eqref{2-7} and the definition of $\gamma_X$ after \eqref{2-14} we see that
\[
Y \geq c^\delta \frac{(1+|s|)\gamma_X}{\gamma} \geq c^\delta(1+|s|) \geq d|s|+B = \xi
\]
if $c^\delta$ is large enough, hence
\[
\begin{split}
\int_{V^\delta}^\infty e^{-\gamma_X u} u^{d|\si|+B} \d u &= \frac{1}{\gamma_X^{d|\si|+B+1}} \int_{V^\delta\gamma_X}^\infty e^{-u} u^{d|\si|+B} \d u \\
&\ll \frac{e^{-V^\delta\gamma_X}(V^\delta\gamma_X)^{d|\si|+B}}{\gamma_X^{d|\si|+B+1}} (1+|\si|) = e^{-V^\delta\gamma_X}V^{\delta(d|\si|+B)} \frac{1+|\si|}{\gamma_X}.
\end{split}
\]
Therefore, recalling that $\delta'<1$, \eqref{2-15} becomes 
\[
\begin{split}
I^{(2)}_X(s,y) &\ll A^{|\si|} e^{C(1+|s|)} \left(e^{-V^\delta\gamma_X}V^{\delta(d|\si|+B)} \frac{1+|\si|}{\gamma_X} + V^{d|\si|+B} \frac{e^{-V^\delta\gamma_X}}{\gamma_X}\right) \\
&\ll \frac{1}{\gamma_X} e^{C(1+|s|)-V^\delta\gamma_X} V^{d|\si|+B}.
\end{split}
\]
But, as before, if $c^\delta$ is large enough we have
\[
\frac12 V^\delta\gamma_X \geq \frac12 c^\delta(1+|s|) \geq C(1+|s|)\quad \text{and} \quad \frac12 V^\delta \geq \frac{c^\delta}{2\gamma} (1+|s|) \geq K,
\]
hence, recalling the definition of $\gamma_X$ after \eqref{2-14} and that $\frac{\beta^{1/d}}{d}y|\omega_X|<1$ for $X>c/\gamma$, we obtain
\begin{equation}
\label{2-17}
I^{(2)}_X(s,y) \ll \frac{1}{\gamma_X} V^{d|\si|+B} y^K.
\end{equation}

\medskip
Gathering \eqref{2-13}, \eqref{2-14} and \eqref{2-17} we finally obtain, in the case under observation, that uniformly for $X>c/\gamma$ 
\begin{equation}
\label{2-18}
I_X(s,y) \ll A^{|\sigma|} e^{C(|s|^{\delta'}+|s|V^{\delta-1})} V^{d|\si|+B} \frac{y^K}{\gamma},
\end{equation}
since $1/\gamma_X \ll 1/\gamma$ for  $X>c/\gamma$ (and the constant in the $\ll$-symbol has the same features as the constant $C$).

\bigskip
{\bf Case $y>d/\beta^{1/d}$.} In this case we change the path of integration $\LL_\infty$ to the half-line $\LL''_\infty$ where $w=u+iV$, $-K\leq u<\infty$, and arguing as in the previous case we obtain
\[
I_X(s,y) = \frac{1}{2\pi i}\int_{\LL''_\infty} h(w,s) (y\omega_X)^{-w}\d w.
\]
Moreover, writing $d_0=d\max_{1\leq j\leq r}|\mu_j|/\lambda_j$ we have for $-K\leq u<-d\si-d_0$ (recall that $\si\leq3/2$) that
\[
\Re(\lambda_js+\mu_j+\frac{\lambda_jw}{d}) <0 \quad \text{and} \quad \Re(\lambda_j(1-s)+\overline{\mu}_j- \frac{\lambda_jw}{d})> 0,
\]
and we split the integral over $\LL''_\infty$ as
\begin{equation}
\label{2-19}
\begin{split}
I_X(s,y) &=\frac{1}{2\pi i}\left(\int_{-K+iV}^{-d\si-d_0+iV} + \int_{-d\si-d_0+iV}^{V^\delta+iV} +  \int_{V^\delta+iV}^{\infty+iV}\right) h(w,s) (y\omega_X)^{-w} \d w \\
&= I^{(3)}_X(s,y) + I^{(4)}_X(s,y) + I^{(5)}_X(s,y),
\end{split}
\end{equation}
say.

\medskip
The treatment of $I^{(3)}_X(s,y)$ is similar to the treatment of the previous case, in the sense that we start from \eqref{2-6} and use \eqref{2-5} to estimate the first two $\Gamma$-factors in the product (with the choice $z=-\frac{\lambda_jw}{d}$ and $a=\lambda_j(1-s)+\overline{\mu}_j$ or $a=1-\lambda_js-\mu_j$, respectively). A computation shows that, assuming \eqref{2-7}, for $w\in\LL_\infty''$ the conditions in \eqref{2-4} are satisfied, and
\[
\begin{split}
\Gamma(\lambda_j(1-s)+\overline{\mu}_j-\frac{\lambda_jw}{d}) &\ll e^{\lambda_j(\frac{V}{d}+t)\arg(-w) +\frac{\lambda_ju}{d}} \big|\frac{\lambda_jw}{d}\big|^{-\lambda_j(\si+\frac{u}{d})+B} e^{C(1+|s|^{\delta'})} \\
\Gamma(1-\lambda_js-\mu_j-\frac{\lambda_jw}{d}) &\ll e^{\lambda_j(\frac{V}{d}+t)\arg(-w) +\frac{\lambda_ju}{d}} \big|\frac{\lambda_jw}{d}\big|^{-\lambda_j(\si+\frac{u}{d})+B} e^{C(1+|s|^{\delta'})}.
\end{split}
\]
Moreover
\[
S(w,s) \ll e^{\frac{\pi}{2}d(\frac{V}{d}+t)}
\]
and by Stirling's formula
\[
|\Gamma(w)| \ll e^{-V\arg w -u} |w|^{u-1/2}.
\]
Therefore, similarly as in \eqref{2-12} and observing that $\arg(-w)-\arg(w) = -\pi$, we obtain
\[
h(w,s) \ll A^{|\si|} e^{C(1+|s|^{\delta'})} e^{dt(\frac{\pi}{2} + \arg(-w))} e^{-\frac{\pi}{2}V} |w|^{-d\si +B} (\frac{d}{\beta^{1/d}})^u,
\]
and hence, since $|(y\omega_X)^{-w}| \leq (y|\omega_X|)^{-u} e^{\frac{\pi}{2}V}$, we get
\[
I^{(3)}_X(s,y) \ll A^{|\si|} e^{C(1+|s|^{\delta'})} \int_{-K}^{-d\si-d_0} e^{dt(\frac{\pi}{2} + \arg(-w))} |w|^{-d\si +B} (\frac{d}{\beta^{1/d}y|\omega_X|})^u \d u.
\]
But, thanks to \eqref{2-7}, $|w|\ll V$ and $\arg(-w)= -\frac{\pi}{2} +O(V^{\delta-1})$. Thus, recalling the definition of $\gamma_X$ after \eqref{2-14} and that $y>d/\beta^{1/d}$ in this case (and hence $\frac{d}{\beta^{1/d}y|\omega_X|}<1$), we have
\begin{equation}
\label{2-20}
\begin{split}
I^{(3)}_X(s,y) &\ll A^{|\si|} e^{C(1+|s|^{\delta'}+|s|V^{\delta-1})} V^{d|\si|+B} \int_{-K}^\infty e^{-\gamma_X u}\d u \\
&\ll A^{|\si|} e^{C(1+|s|^{\delta'}+|s|V^{\delta-1})} V^{d|\si|+B} \frac{y^K}{\gamma_X}.
\end{split}
\end{equation}

\medskip
In order to treat $I^{(4)}_X(s,y)$ we use the reflection formula of the $\Gamma$ function to write
\[
h(w,s) = \pi^r \Gamma(w) \widetilde{S}(w,s) \prod_{j=1}^r \frac{1}{\Gamma(\lambda_js +\mu_j +\frac{\lambda_jw}{d}) \Gamma(1-\lambda_j(1-s)-\overline{\mu}_j + \frac{\lambda_jw}{d})}
\]
with
\[
\widetilde{S}(w,s) = \prod_{j=1}^r \frac{1}{\sin\pi(\lambda_j(1-s-\frac{w}{d})+\overline{\mu}_j)}.
\]
By the factorial formula of the $\Gamma$ function we rewrite $h(w,s)$ as
\begin{equation}
\label{2-21}
h(w,s) = \pi^r \Gamma(w) \widetilde{S}(w,s) P(s+\frac{w}{d})\prod_{j=1}^r \frac{1}{\Gamma(a_j(s) +\frac{\lambda_jw}{d}) \Gamma(b_j(s)+ \frac{\lambda_jw}{d})}
\end{equation}
where $P\in\CC[z]$ has $\deg P = 2r\nu_0$,
\[
a_j(s) = \lambda_js +\mu_j + \nu_0, \hskip1.5cm b_j(s) = 1-\lambda_j(1-s)-\overline{\mu}_j + \nu_0
\]
and $\nu_0=\nu_0(F)\in\NN$ is such that for $u\geq -d\si-d_0$
\[
\Re(a_j(s) +\frac{\lambda_jw}{d}), \Re(b_j(s)+ \frac{\lambda_jw}{d}) \geq 0.
\]
Now we apply \eqref{2-5} with $z=\frac{\lambda_jw}{d}$ and $a=a_j(s), b_j(s)$ respectively. A computation shows that conditions \eqref{2-4} are satisfied thanks to \eqref{2-7} and to the choice of $\nu_0$, and
\[
\begin{split}
\Gamma(a_j(s) + \frac{\lambda_jw}{d}) &\gg e^{-\lambda_j(t+\frac{V}{d})\arg w -\frac{\lambda_ju}{d}} \big|\frac{\lambda_j w}{d}\big|^{\lambda_j(\si+\frac{u}{d}) +B} e^{C(1+|s|^{\delta'})} \\
\Gamma(b_j(s) + \frac{\lambda_jw}{d}) &\gg e^{-\lambda_j(t+\frac{V}{d})\arg w -\frac{\lambda_ju}{d}} \big|\frac{\lambda_j w}{d}\big|^{\lambda_j(\si+\frac{u}{d}) +B} e^{C(1+|s|^{\delta'})}.
\end{split}
\]
Moreover, thanks to \eqref{2-7} and to Stirling's formula we have
\[
\widetilde{S}(w,s) \ll e^{-\frac{\pi}{2}d(\frac{V}{d}+t)}, \hskip.5cm P(s+\frac{w}{d}) \ll |w|^B, \hskip.5cm \Gamma(w) \ll e^{-V\arg w -u} |w|^{u-1/2},
\]
hence from \eqref{2-21} we obtain
\[
h(w,s) \ll A^{|\si|} e^{C(1+|s|^{\delta'})} e^{-\frac{\pi}{2}V +dt(\arg w-\frac{\pi}{2})} |w|^{-d\si +B} (\frac{d}{\beta^{1/d}})^u.
\]
As before we have $|(y\omega_X)^{-w}| \leq (y|\omega_X|)^{-u} e^{\frac{\pi}{2}V}$, thus for $u\geq -d\si-d_0$
\begin{equation}
\label{2-22}
h(w,s) (y\omega_X)^{-w} \ll A^{|\si|} e^{C(1+|s|^{\delta'})} e^{dt(\arg w-\frac{\pi}{2})} |w|^{-d\si +B} (\frac{\beta^{1/d}y|\omega_X|}{d})^{-u}.
\end{equation}
But for $-d\si-d_0 \leq u \leq V^\delta$ we have $\arg w = \frac{\pi}{2} + O(V^{\delta-1})$ and $|w|\ll V$, therefore recalling the definition of $\gamma_X$ after \eqref{2-14} and that $y>d/\beta^{1/d}$ we obtain
\begin{equation}
\label{2-23}
\begin{split}
I^{(4)}_X(s,y) &\ll A^{|\si|} e^{C(1+|s|^{\delta'}+|s|V^{\delta-1})} V^{d|\si|+B} \int_{-K}^\infty e^{-u\gamma_X} \d u \\
&\ll A^{|\si|} e^{C(1+|s|^{\delta'}+|s|V^{\delta-1})} V^{d|\si|+B} \frac{y^K}{\gamma_X}.
\end{split}
\end{equation}

\medskip
Finally, since $\arg w \leq \frac{\pi}{2} + O(V^{\delta-1})$ for $u\geq V^\delta$, from \eqref{2-22} we get
\[
I^{(5)}_X(s,y) \ll A^{|\si|} e^{C(1+|s|^{\delta'}+|s|V^{\delta-1})} \int_{V^\delta}^\infty e^{-u\gamma_X} (u^{d|\si|+B} + V^{d|\si|+B})\d u,
\]
and by the same argument used to estimate $I^{(2)}_X(s,y)$ we obtain
\begin{equation}
\label{2-24}
I^{(5)}_X(s,y) \ll A^{|\si|} e^{C(1+|s|^{\delta'}+|s|V^{\delta-1})} V^{d|\si|+B} \frac{y^K}{\gamma_X}.
\end{equation}

\medskip
From \eqref{2-19}, \eqref{2-20}, \eqref{2-23} and \eqref{2-24} we have, in the case under observation, that uniformly for $X> c/\gamma$ (see after \eqref{2-18})
\begin{equation}
\label{2-25}
I_X(s,y) \ll A^{|\si|} e^{C(|s|^{\delta'}+|s|V^{\delta-1})} V^{d|\si|+B} \frac{y^K}{\gamma},
\end{equation}
and in view of \eqref{2-18} we have that \eqref{2-25} holds for every $y>0$.

\medskip
We are now ready to conclude the proof of Lemma B. From \eqref{2-2}, \eqref{2-3}, \eqref{2-7} and \eqref{2-25} we have, uniformly for $X> c/\gamma$, that
\[
F^{(3)}_X(s) \ll A^{|\si|} e^{C(|s|^{\delta'}+|s|V^{\delta-1})} V^{d|\si|+B} \sum_{n=1}^\infty \frac{|a(n)|}{n^{1-\sigma}n^{K/d}|\log(n_\alpha/n)|}.
\]
Lemma B follows since $n\neq n_\alpha$ and the series is convergent for every $\sigma\leq 3/2$.

\bigskip
\section{Proof of Theorem 2}

\smallskip
In this section we follow the notation of Section 1 and suppose that $F(s),\alpha, \delta_0,\delta,\eta$ are as in Theorem 2. Moreover, we always assume (including in the statement of the lemmas) that $\si\leq 3/2$ and $|s|\leq R$ with an arbitrary $R\geq 1$. The initial steps of the proof of Theorem 2 are identical to those of Theorem 1. Indeed, we borrow Lemma A, whose conclusions remain unchanged under the hypotheses of Theorem 2; see the remark after the proof of Lemma A. Moreover, the conclusions of Lemma B also remain unchanged under the hypotheses of Theorem 2, provided we replace $F_X^{(3)}(s)$ by $\widetilde{F}_X^{(3)}(s)$, where, recalling \eqref{2-2} and the definition of $I_X(s,y)$ before \eqref{2-2},
\[
\widetilde{F}_X^{(3)}(s) = \omega Q^{1-2s} \sum_{n\neq n_\alpha}\frac{\overline{a(n)}}{n^{1-s}} I_X(s, y_n).
\]
Therefore, with the notation of Section 2 we have for $\sigma\leq 3/2$
\begin{equation}
\label{3-1}
F_X(s,\alpha) = F_X^{(1)}(s) + F_X^{(2)}(s) +  \widetilde{F}_X^{(3)}(s) + \omega Q^{1-2s} \frac{\overline{a(n_\alpha)}}{n_\alpha^{1-s}} I_X(s, d/\beta^{1/d})
\end{equation}
and, provided $V\geq V_0(1+|s|)^{1/\delta}$,
\begin{equation}
\label{3-2}
F_X^{(1)}(s) + F_X^{(2)}(s) +  \widetilde{F}_X^{(3)}(s) \ll A^{|\si|} V^{d|\si|+B} e^{C(|s|^{\delta'}+|s|V^{\delta-1})}
\end{equation}
uniformly for $X\geq X_0$. The rest of the proof is devoted to the evaluation and estimation of $I_X(s, d/\beta^{1/d})$, which for convenience we denote simply by $I_X(s)$. 

\medskip
For $w\in\LL_\infty$, as in \eqref{2-6} we write
\[
h(w,s) = S(s+\frac{w}{d}) \Gamma(w) \prod_{j=1}^r \Gamma(a_j-\lambda_j\frac{w}{d}) \Gamma(b_j-\lambda_j\frac{w}{d})
\]
with
\[
S(z) = \pi^{-r} \prod_{j=1}^r \sin \pi(\lambda_jz+\mu_j), \hskip1cm a_j=\lambda_j(1-s) + \overline{\mu}_j, \hskip1cm b_j=1-\mu_j-\lambda_js.
\]
We first note that a standard calculation based on the expression $\sin z = \frac{e^{iz} - e^{-iz}}{2i}$ gives for $w\in\LL_\infty$ (recall that $v\geq V\geq |s|$)
\begin{equation}
\label{3-3}
S(s+\frac{w}{d}) = \frac{1}{(2\pi)^r} e^{-i\frac{\pi}{2}\xi_F} e^{-id\frac{\pi}{2}(s+\frac{w}{d})} \big(1+ O(e^{-\lambda_0v})\big)
\end{equation}
with some $0<\lambda_0\leq1$. We recall that the $H$-invariants of a function $F\in\S^\sharp$ are defined as
\[
H_F(n) = 2\sum_{j=1}^r\frac{B_n(\mu_j)}{\lambda_j^{n-1}} \hskip2cm n=0,1,...
\]
where $B_n(x)$ is the $n$-th Bernoulli polynomial; see \cite{Ka-Pe/2002a} for the properties of such invariants. Denoting by $R_\nu(s)$ the polynomials
\begin{equation}
\label{3-4}
\begin{split}
R_\nu(s) &= 2^\nu\big(B_{\nu+1}(a+\frac12) - B_{\nu+1}(\frac12-b)\big) \\
&+ \frac{d^\nu}{2}  \sum_{k=0}^{\nu+1}{\nu+1\choose k}\big((-1)^\nu H_F(k)s^{\nu+1-k} - \overline{H_F(k)}(1-s)^{\nu+1-k}\big),
\end{split}
\end{equation}
where $a=a(s) = \frac{d}{2}(1-s)+\frac12\overline{\xi_F}$ and $b=b(s) = \frac{d}{2}s+\frac12\xi_F$, we have

\medskip
{\bf Lemma 3.1.} {\sl Let $w\in\LL_\infty$, $|s|\leq R$, $\si\leq 3/2$, $1\leq M\leq DR$ and $V\geq (c_1R)^{1/\delta}$, where $D>0$ is arbitrary and $c_1$ may depend also on $D$. Then we have
\[
h(w,s) = c_0(F)\frac{\Gamma(a+\frac12 -\frac{w}{2})}{\Gamma(b+\frac12+\frac{w}{2})} \Gamma(w) \big(\frac{d}{2\beta^{1/d}}\big)^{ds+w} \exp\big(\sum_{\nu=1}^M \frac{R_\nu(s)}{\nu(\nu+1)} \frac{1}{w^\nu} + \rho_M\big),
\]
with
\[
\rho_M \ll c_2^M\frac{R^{M+2}}{|w|^{M+1}} + e^{-\lambda_0v},
\]
where $c_0(F)\neq 0$, $\lambda_0>0$ is a certain constant, $c_2$ may depend also on $D$, and we may assume that $c_1\geq c_2$.}

\medskip
{\it Proof.} We use Lemma D of Section 2 with $z=-\lambda_j\frac{w}{d}$, $a=a_j$ or $a=b_j$ and $N=M$ to get an asymptotic expansion of $\Gamma(a_j-\lambda_j\frac{w}{d}) \Gamma(b_j-\lambda_j\frac{w}{d})$. Since $u=-K$ and $\si\leq 3/2$, and thanks to the restrictions on $V$ and $M$, for $w\in\LL_\infty$ the hypotheses of Lemma D are satisfied, hence
\[
\begin{split}
\log\Gamma(a_j-\lambda_j\frac{w}{d}) \Gamma(b_j-\lambda_j\frac{w}{d}) &= (\lambda_j-2\lambda_j(s+\frac{w}{d}) -2i\Im\mu_j)\log(-\frac{w}{2d}) \\
&+ (\lambda_j-2\lambda_j(s+\frac{w}{d}) -2i\Im\mu_j)\log(2\lambda_j) +2\lambda_j\frac{w}{d} + \log2\pi \\
&-\sum_{\nu=1}^M d^\nu \frac{B_{\nu+1}(a_j)+B_{\nu+1}(b_j)}{\nu(\nu+1)} \frac{1}{\lambda_j^\nu w^\nu} + O(\rho_{M,j})
\end{split}
\]
where, thanks to the restriction on $V$,
\[
\rho_{M,j} \ll \left(\frac{d}{\lambda_j|w|}\right)^{M+1}\left((M+\frac{(\lambda_j|s|+B)^2}{M^2})(\lambda_j|s|+B)^M+B^MM!\right) \ll c^M\frac{R^{M+2}}{|w|^{M+1}}.
\]
Therefore, summing over $j=1,\dots,r$ we obtain
\begin{equation}
\label{3-5}
\begin{split}
\log\prod_{j=1}^r &\left(\Gamma(a_j-\lambda_j\frac{w}{d}) \Gamma(b_j-\lambda_j\frac{w}{d})\right) = (\frac{d}{2} -ds-w -id\theta_F)\log(-\frac{w}{2}) +c_1(F) \\
&+(ds+w)\log(\frac{d}{2\beta^{1/d}}) +w -\sum_{\nu=1}^M \frac{\widetilde{R}_\nu(s)}{\nu(\nu+1)}\frac{1}{w^\nu} + O(c^M\frac{R^{M+2}}{|w|^{M+1}}),
\end{split}
\end{equation}
with a certain $c_1(F)$ and
\[
\widetilde{R}_\nu(s) = d^\nu\sum_{j=1}^r \frac{B_{\nu+1}(a_j)+B_{\nu+1}(b_j)}{\lambda_j^\nu}.
\]
Now we note that
\[
\frac{\Gamma(a+\frac12 -\frac{w}{2})}{\Gamma(b+\frac12+\frac{w}{2})} = \frac{\sin(\pi(b+\frac12+\frac{w}{2}))}{\pi} \Gamma(a+\frac12 -\frac{w}{2}) \Gamma(\frac12 - b - \frac{w}{2})
\]
and, since $w\in\LL_\infty$ and hence $\Im(b+\frac12+\frac{w}{2})>0$,
\[
\sin(\pi(b+\frac12+\frac{w}{2})) = -\frac{1}{2i} e^{-i\pi(b+\frac12+\frac{w}{2})} (1+O(e^{-v})) = \frac12e^{-i\frac{\pi}{2}\xi_F} e^{-i\frac{\pi d}{2}(s+\frac{w}{d})} (1+O(e^{-v})).
\]
Moreover, arguing as before, from Lemma D we get
\[
\begin{split}
\log \Gamma(a+\frac12 -\frac{w}{2}) &\Gamma(\frac12 - b - \frac{w}{2}) = (\frac{d}{2} -ds-w-id\theta_F)\log(-\frac{w}{2}) +w + \log(2\pi) \\
&- \sum_{\nu=1}^M \frac{2^\nu(B_{\nu+1}(a+\frac12)-B_{\nu+1}(\frac12-b))}{\nu(\nu+1)} \frac{1}{w^\nu} + O(c^M\frac{R^{M+2}}{|w|^{M+1}}).
\end{split}
\]
Hence, with a certain $c_2(F)$,
\begin{equation}
\label{3-6}
\begin{split}
\log \frac{\Gamma(a+\frac12 -\frac{w}{2})}{\Gamma(b+\frac12+\frac{w}{2})} &= -i\frac{\pi d}{2}(s+\frac{w}{d}) + (\frac{d}{2} -ds-w-id\theta_F)\log(-\frac{w}{2}) +w +c_2(F) \\
&- \sum_{\nu=1}^M \frac{2^\nu(B_{\nu+1}(a+\frac12)-B_{\nu+1}(\frac12-b))}{\nu(\nu+1)} \frac{1}{w^\nu} + O(c^M\frac{R^{M+2}}{|w|^{M+1}}+e^{-v}).
\end{split}
\end{equation}
Recalling the definition of $h(w,s)$ and \eqref{3-3}, comparing \eqref{3-5} and \eqref{3-6} we obtain
\begin{equation}
\label{3-7}
\begin{split}
\log h(w,s) &= c_3(F) \log \frac{\Gamma(a+\frac12 -\frac{w}{2})}{\Gamma(b+\frac12+\frac{w}{2})} \log\Gamma(w) (ds+w)\log(\frac{d}{2\beta^{1/d}}) \\
&+ \sum_{\nu=1}^M \frac{R^*_\nu(s)}{\nu(\nu+1)} \frac{1}{w^\nu} + O(c^M\frac{R^{M+2}}{|w|^{M+1}}+e^{-\lambda_0v})
\end{split}
\end{equation}
with a certain $c_3(F)$ and
\begin{equation}
\label{3-8}
R^*_\nu(s) = 2^\nu(B_{\nu+1}(a+\frac12)-B_{\nu+1}(\frac12-b)) - d^\nu\sum_{j=1}^r \frac{B_{\nu+1}(a_j)+B_{\nu+1}(b_j)}{\lambda_j^\nu};
\end{equation}
Lemma 3.1 follows from \eqref{3-7} as soon as we prove that $R^*_\nu(s)=R_\nu(s)$. But from the properties of the Bernoulli polynomials, see Section 1.13 of Bateman's Project \cite{EMOT/1953}, we have
\[
\begin{split}
\sum_{j=1}^r &\frac{B_{\nu+1}(a_j)+B_{\nu+1}(b_j)}{\lambda_j^\nu} = \sum_{j=1}^r \frac{B_{\nu+1}(\lambda_j(1-s) + \overline{\mu}_j)+(-1)^{\nu+1}B_{\nu+1}(\lambda_js+\mu_j)}{\lambda_j^\nu} \\
&= \sum_{j=1}^r \sum_{k=0}^{\nu+1} \frac{{\nu+1\choose k}\lambda_j^{\nu+1-k} B_{\nu+1}(\overline{\mu}_j)(1-s)^{\nu+1-k} +(-1)^{\nu+1}{\nu+1\choose k}\lambda_j^{\nu+1-k} B_{\nu+1}(\mu_j)s^{\nu+1-k}}{\lambda_j^\nu} \\
&= \sum_{k=0}^{\nu+1} {\nu+1 \choose k} \left(\frac12\overline{H_F(k)}(1-s)^{\nu+1-k} +(-1)^{\nu+1}\frac12H_F(k)s^{\nu+1-k}\right) \\
&=-\frac12 \sum_{k=0}^{\nu+1} {\nu+1 \choose k} \left((-1)^\nu H_F(k)s^{\nu+1-k} -\overline{H_F(k)}(1-s)^{\nu+1-k} \right),
\end{split}
\]
hence $R^*_\nu(s)=R_\nu(s)$ in view of \eqref{3-4} and \eqref{3-8}; see also Lemma 3.3 of \cite{Ka-Pe/twist}. \fine

\medskip
The treatment of $\exp\big(\sum_{\nu=1}^M \frac{R_\nu(s)}{\nu(\nu+1)} \frac{1}{w^\nu} + \rho_M\big)$ in Lemma 3.1 is partly similar to the arguments in Section 3 of \cite{Ka-Pe/twist}, but the differences are such that we cannot simply quote the results in \cite{Ka-Pe/twist}. However, we will be a bit more sketchy in such a treatment, and we shall refer to \cite{Ka-Pe/twist} whenever possible.

\medskip
{\bf Lemma 3.2.} {\sl Let $R_\nu(s)$ be as in \eqref{3-4} and $1\leq \nu \leq DR$, where $D>0$ is arbitrary. Then
\[
R_\nu(s) \ll (c_3R)^{\nu+1}
\]
where $c_3$ may depend also on $D$.}

\medskip
{\it Proof.} For $d=2$ this is essentially Lemma 3.8 of \cite{Ka-Pe/twist} (notice the slightly different definition of $R_\nu(s)$ in \cite{Ka-Pe/twist}), and the proof in the general case is similar. \fine

\medskip
{\bf Lemma 3.3.} {\sl Under the hypotheses and with the notation of Lemma $3.1$, and assuming in addition that $\frac{2\delta-1}{1-\delta} \leq M$, we have
\[
\begin{split}
h(w,s) &= c_0(F)\frac{\Gamma(a+\frac12 -\frac{w}{2})}{\Gamma(b+\frac12+\frac{w}{2})} \Gamma(w) \big(\frac{d}{2\beta^{1/d}}\big)^{ds+w} \exp\big(\sum_{\nu=1}^M \frac{R_\nu(s)}{\nu(\nu+1)} \frac{1}{w^\nu}\big) \\
& + O\big(A^{|\si|}e^{O(R^2/V)} \frac{c_4^MR^{M+2}|\Gamma(w)|}{|w|^{M+d\si-K-d/2+1}}\big),
\end{split}
\]
where $A$ and $c_4$ may depend also on $D$.}

\medskip
{\it Proof.} By Lemma 3.1 we have
\[
\rho_M \ll \frac{(c_2R)^{M+2}}{(c_1R)^{(M+1)/\delta}} + e^{-\lambda_0v} \ll R^{M+2-(M+1)/\delta} + e^{-\lambda_0v},
\]
and thanks to the restriction $\frac{2\delta-1}{1-\delta} \leq M$ we see that the exponent of $R$ is negative. Therefore $\rho_M\ll1$ and hence $\exp(\rho_M) = 1 +O(\rho_M)$. Moreover, for $w\in\LL_\infty$ we have
\[
e^{-\lambda_0v} \ll e^{-\lambda_0|w|} \ll \frac{(M+1)!}{(\lambda_0|w|)^{M+1}} \ll \left(\frac{M+1}{\lambda_0|w|}\right)^{M+1} \ll \frac{(cR)^{M+2}}{|w|^{M+1}}.
\]
Hence in order to prove the lemma we have to show that
\begin{equation}
\label{3-9}
\left|\frac{\Gamma(a+\frac12 -\frac{w}{2})}{\Gamma(b+\frac12+\frac{w}{2})} \big(\frac{d}{2\beta^{1/d}}\big)^{ds+w} \exp\big(\sum_{\nu=1}^M \frac{R_\nu(s)}{\nu(\nu+1)} \frac{1}{w^\nu}\big) \right| \ll A^{|\si|}\frac{e^{O(R^2/V)}}{|w|^{d\si-K-d/2}}.
\end{equation}
But from \eqref{3-6} with $M=1$ we obtain
\begin{equation}
\label{3-10}
\left|\frac{\Gamma(a+\frac12 -\frac{w}{2})}{\Gamma(b+\frac12+\frac{w}{2})}\right| \ll A^{|\si|} |w|^{\frac{d}{2} -d\si +K} e^{O(R^2/V)},
\end{equation}
$\big(\frac{d}{2\beta^{1/d}}\big)^{ds+w} \ll A^{|\si|}$ and from Lemma 3.2 (observing that we may choose $c_1\geq c_3$) we have
\[
\sum_{\nu=1}^M \frac{|R_\nu(s)|}{\nu(\nu+1)} \frac{1}{|w|^\nu} \ll R\sum_{\nu=1}^M \frac{1}{\nu(\nu+1)}(\frac{c_3R}{v})^\nu \ll \frac{R^2}{V}+1.
\]
Therefore \eqref{3-9} follows and the proof of Lemma 3.3 is complete. \fine

\medskip
With the notation of Lemma 3.1 and writing
\begin{equation}
\label{3-11}
\begin{split}
h_M^*(w,s) &= \frac{\Gamma(a+\frac12 -\frac{w}{2})}{\Gamma(b+\frac12+\frac{w}{2})} \Gamma(w) \exp\big(\sum_{\nu=1}^M \frac{R_\nu(s)}{\nu(\nu+1)} \frac{1}{w^\nu}\big) \\
I^*_{X.M}(s) &= \frac{1}{2\pi i} \int_{\LL_\infty} h_M^*(w,s) (2\omega_X)^{-w} \d w,
\end{split}
\end{equation}
from Lemma 3.3 we get
\begin{equation}
\label{3-12}
I_X(s) = c_0(F) \big(\frac{d}{2\beta^{1/d}}\big)^{ds} I^*_{X.M}(s) + O(A^{|\si|} e^{O(R^2/V)} (cR)^{M+2})
\end{equation}
uniformly for $X\geq X_0$, provided $V\geq (c_1R)^{1/\delta}$ and $\max(\frac{2\delta-1}{1-\delta}, d|\si|+d/2) \leq M \leq DR$. Indeed, for $w\in\LL_\infty$ we have
\[
\Gamma(w) \ll e^{-\frac{\pi}{2}v} |w|^{-K-1/2} \qquad \text{and} \qquad \big|(d\omega_X/\beta^{1/d})^{-w}\big| \ll e^{\frac{\pi}{2}v},
\]
hence replacing $h(w,s)$ in $I_X(s)$ by its main term in Lemma 3.3 we get an error of size
\[
\ll A^{|\si|} e^{O(R^2/V)} c^M R^{M+2}   \int_V^\infty \frac{\d v}{v^{M+d\si-d/2+3/2}} \ll A^{|\si|} e^{O(R^2/V)} (cR)^{M+2}
\]
since $M\geq  d|\si|+d/2$, and \eqref{3-12} follows.

\medskip
In order to study $h^*_M(w,s)$ we write $V_0(s)=1$ identically and for $\mu\geq1$
\begin{equation}
\label{3-13}
V_\mu(s) = \sum_{1\leq m\leq \mu} \frac{1}{m!} \sum_{\substack{\nu_1\geq1,...,\nu_m\geq1\\ \nu_1+...+\nu_m=\mu}} \prod_{j=1}^m\frac{R_{\nu_j}(s)}{\nu_j(\nu_j+1)}.
\end{equation}

\medskip
{\bf Lemma 3.4.} {\sl Under the hypotheses and with the notation of Lemma $3.1$ we have
\[
\exp\big(\sum_{\nu=1}^M \frac{R_\nu(s)}{\nu(\nu+1)} \frac{1}{w^\nu}\big) = \sum_{\mu=0}^M \frac{V_\mu(s)}{w^\mu} + O\big(\frac{(c_5R)^M}{|w|^{\delta M}} e^{O(RV^{\delta-1})}\big),
\]
where $c_5$ may depend also on $D$.}

\medskip
{\it Proof.} Expanding the exponential and recalling \eqref{3-13} we have
\[
\exp\big(\sum_{\nu=1}^M \frac{R_\nu(s)}{\nu(\nu+1)} \frac{1}{w^\nu}\big) = \sum_{\mu=0}^M \frac{V_\mu(s)}{w^\mu} + E 
\]
with
\[
E = \sum_{\mu>M} \frac{V_{\mu,M}(s)}{w^\mu}, \hskip1.5cm V_{\mu,M}(s) = \sum_{1\leq m\leq \mu} \frac{1}{m!} \sum_{\substack{1\leq \nu_j\leq M\\ \nu_1+...+\nu_m=\mu}} \prod_{j=1}^m\frac{R_{\nu_j}(s)}{\nu_j(\nu_j+1)}.
\]
By Lemma 3.2 we obtain
\[
V_{\mu,M}(s) \ll \sum_{1\leq m\leq \mu} \frac{1}{m!} \sum_{\substack{1\leq \nu_j\leq M\\ \nu_1+...+\nu_m=\mu}} (cR)^{\mu+m} \prod_{j=1}^m\frac{1}{\nu_j(\nu_j+1)} \ll (cR)^\mu \sum_{m=0}^\mu \frac{(cR)^m}{m!},
\]
hence for $\mu\leq \frac{3}{2}cR$ we have
\begin{equation}
\label{3-14}
V_{\mu,M}(s) \ll \frac{(cR)^{2\mu}}{\mu!}
\end{equation}
thanks to Lemma 3.6 of \cite{Ka-Pe/twist}, while for $\mu> \frac{3}{2}cR$ we get
\[
V_{\mu,M}(s) \ll (cR)^\mu e^{cR}.
\]
Therefore we have
\[
E \ll \sum_{M<\mu\leq 3cR/2} \frac{(cR)^{2\mu}}{|w|^\mu\mu!} + \sum_{\mu>3cR/2} \frac{(cR)^\mu e^{cR}}{|w|^\mu} = \Sigma_1+\Sigma_2,
\]
say. But, thanks to the restrictions on $w$ and $V$,
\[
\Sigma_1 \ll \sum_{\mu>M} \frac{(cR)^\mu}{|w|^{\delta\mu}} \frac{(cR)^\mu}{\mu!V^{(1-\delta)\mu}} \ll \frac{(cR)^M}{|w|^{\delta M}} e^{O(RV^{\delta-1})}
\]
and
\[
\Sigma_2 \ll  e^{cR} \sum_{\mu>3cR/2} \left(\frac{cR}{|w|}\right)^\mu \ll e^{cR} \left(\frac{cR}{|w|}\right)^{3cR/2} \ll \left(\frac{cR}{|w|}\right)^M
\]
since $M<\frac{3}{2}cR$. Lemma 3.4 is therefore proved. \fine

\medskip
Now we transform the sum over $\mu$ in Lemma 3.4 by arguments similar to those in Section 3 of \cite{Ka-Pe/twist}, see Lemmas 3.11-3.18. We start with the following variant of Lemma 3.13 of \cite{Ka-Pe/twist}, asserting that for $|w|\geq 2M$ and $1\leq \mu\leq M$
\begin{equation}
\label{3-15}
\frac{1}{w^\mu} = \sum_{\ell =\mu}^M \frac{C_{\mu,\ell}}{(w-1)\cdots(w-\ell)} + O\big(\frac{4^MM!}{(\mu-1)!} \frac{1}{|w|^{M+1}}\big),
\end{equation}
where the coefficients $C_{\mu,\ell}$ satisfy
\begin{equation}
\label{3-16}
|C_{\mu,\ell}| \leq \frac{(\ell-1)!}{(\mu-1)!}{\ell-1\choose \mu-1}
\end{equation}
by Lemma 3.12 of \cite{Ka-Pe/twist}. Indeed, Lemma 3.13 of \cite{Ka-Pe/twist} gives \eqref{3-15} with an error term
\[
\ll \frac{2^MM!}{(\mu-1)!} \frac{1}{|w(w-1)\cdots(w-M)|},
\]
and $|w(w-1)\cdots(w-M)| \geq (|w|-M)^{M+1} \geq (|w|/2)^{M+1}$ thus giving \eqref{3-15}. Let
\begin{equation}
\label{3-17}
A_{\mu,\ell}(s) = \sum_{k=0}^{\ell-\mu} (-1)^{\mu+\ell+k} C_{\mu+k,\ell} {-\mu\choose k} (s+1)^k.
\end{equation}

\medskip
{\bf Lemma 3.5.} {\sl Let $1\leq M \leq DR$ with arbitrary $D\geq 1$, $1\leq \mu\leq M$ and $|w|\geq \max(2M,3R+4)$. Then with the above notation}
\[
\frac{1}{(w+s)^\mu} = \sum_{\ell=\mu}^M \frac{A_{\mu,\ell}(s)}{w(w+1)\cdots (w+\ell-1)} +O\big(24^M\frac{(DR)^{M-\mu+1}}{|w|^{M+1}}\big).
\]

\medskip
{\it Proof.} From \eqref{3-15} we have
\[
\frac{(-1)^\mu}{w^\mu} = \frac{1}{(-w)^\mu} = \sum_{\ell =\mu}^M \frac{(-1)^\ell C_{\mu,\ell}}{(w+1)\cdots(w+\ell)} + O\big(\frac{4^MM!}{(\mu-1)!} \frac{1}{|w|^{M+1}}\big),
\]
hence
\begin{equation}
\label{3-18}
\frac{1}{(w-1)^\mu} = \sum_{\ell =\mu}^M \frac{(-1)^{\ell+\mu} C_{\mu,\ell}}{w\cdots(w+\ell-1)} + O\big(\frac{8^MM!}{(\mu-1)!} \frac{1}{|w|^{M+1}}\big).
\end{equation}
Since $|\frac{s+1}{w-1}| \leq \frac{R+1}{|w|-1}\leq \frac{1}{3}$ thanks to the hypothesis $|w|\geq 3R+4$, arguing as in the proof of Lemma 3.14 of \cite{Ka-Pe/twist} and starting with (3.16) there, we have
\[
\frac{1}{(w+s)^\mu} = \frac{1}{(w-1)^\mu(1+\frac{s+1}{w-1})^\mu} = \sum_{k=0}^{M-\mu} {-\mu\choose k} \frac{(s+1)^k}{(w-1)^{\mu+k}} +O(4^M\frac{(R+1)^{M-\mu+1}}{|w|^{M+1}}).
\]
Inserting \eqref{3-18} in the last equation and recalling \eqref{3-17} we obtain
\[
\begin{split}
\frac{1}{(w+s)^\mu} = \sum_{\ell=\mu}^M \frac{A_{\mu,\ell}(s)}{w(w+1)\cdots (w+\ell-1)} &+ O\big(\sum_{k=0}^{M-\mu} \left| {-\mu\choose k}\right| (R+1)^k \frac{8^M M!}{(k+\mu-1)!} \frac{1}{|w|^{M+1}}\big) \\
&+O\big(4^M\frac{(R+1)^{M-\mu+1}}{|w|^{M+1}}\big);
\end{split}
\]
we denote by $E_1$ and $E_2$ the two error terms in the last equation, respectively. Since $|{-\mu\choose k}| = {\mu+k-1\choose k}$, using Lemma 3.6 of \cite{Ka-Pe/twist} we have
\[
\begin{split}
E_1 &\ll \frac{6^MM!}{|w|^{M+1}} \sum_{k=0}^{M-\mu}\frac{(k+\mu-1)!}{k!(\mu-1)!} \frac{(R+1)^k}{(k+\mu-1)!} \ll \frac{9^MM!}{|w|^{M+1}(\mu-1)!} \sum_{k=0}^{M-\mu} \frac{(DR)^k}{k!} \\
&\ll \frac{12^MM!}{|w|^{M+1}}\frac{(DR)^{M-\mu}}{\mu!(M-\mu)!} \ll 24^M \frac{(DR)^{M-\mu}}{|w|^{M+1}}.
\end{split}
\]
The lemma follows recalling the bound for $E_2$. \fine

\medskip
Under the hypotheses of Lemma 3.5, but with  $|w|\geq \max(2M+R,4(R+1))$, we also have
\begin{equation}
\label{3-19}
\frac{1}{w^\mu} = \sum_{\ell=\mu}^M \frac{A_{\mu,\ell}(-s)}{(w+s)\cdots (w+s+\ell-1)} +O\big(50^M\frac{(DR)^{M-\mu+1}}{|w|^{M+1}}\big).
\end{equation}
Indeed, we apply Lemma 3.5 with $w+s$ in place of $w$ and $-s$ in place of $s$, thus getting the main term in \eqref{3-19}  plus an error satisfying (since $|w+s|\geq |w|-|s|\geq |w|/2$)
\[
\ll 25^M\frac{(DR)^{M-\mu+1}}{|w+s|^{M+1}} \ll 50^M\frac{(DR)^{M-\mu+1}}{|w|^{M+1}}.
\]
Recalling  \eqref{3-4}, \eqref{3-17} and the definition of $b=b(s)$ before Lemma 3.1, writing for $\nu\geq 0$
\begin{equation}
\label{3-20}
Q_\nu(s) = \sum_{\mu=0}^\nu \frac{1}{2^\mu} V_\mu(s) A_{\mu,\nu}(-b-\frac12)
\end{equation}
and observing that $Q_0(s)=1$ identically, we have

\medskip
{\bf Lemma 3.6.} {\sl Under the hypotheses and with the notation of Lemma $3.1$ we have
\[
\begin{split}
\exp\big(\sum_{\nu=1}^M \frac{R_\nu(s)}{\nu(\nu+1)} \frac{1}{w^\nu}\big) = 1 + &\sum_{\nu=1}^M \frac{Q_\nu(s)}{(b+\frac12+\frac{w}{2})\cdots (b+\frac12+\frac{w}{2}+\nu-1)} \\
&+O\big(\frac{(c_6R)^M}{|w|^{\delta M}}e^{O(RV^{\delta-1})}\big),
\end{split}
\]
where $c_6$ may depend also on $D$.}

\medskip
{\it Proof.} From \eqref{3-19} applied to $\frac{1}{w^\mu} = \frac{1}{2^\mu} \frac{1}{(\frac{w}{2})^\mu}$ (changing $\ell$ into $\nu$) and recalling \eqref{3-20} we get
\[
\begin{split}
\sum_{\mu=0}^M \frac{V_\mu(s)}{w^\mu}  &= \sum_{\mu=0}^M \frac{V_\mu(s)}{2^\mu}
\sum_{\nu=\mu}^M \frac{A_{\mu,\nu}(-b-\frac12)}{(\frac{w}{2}+b+\frac12)\cdots (\frac{w}{2}+b+\frac12+\nu-1)} \\
&\hskip2cm+O\big(c^M\sum_{\mu=0}^M \frac{|V_\mu(s)|}{2^\mu} \frac{(DR)^{M-\mu+1}}{|w|^{M+1}}\big) \\
&= 1 + \sum_{\nu=1}^M \frac{Q_\nu(s)}{(b+\frac12+\frac{w}{2})\cdots (b+\frac12+\frac{w}{2}+\nu-1)} + E,
\end{split}
\]
say. In view of the restriction on $M$, for $1\leq \mu\leq M$ from \eqref{3-14} we have 
\[
V_\mu(s) \ll \frac{(cR)^{2\mu}}{\mu!}
\]
and hence, since $w\in\LL_\infty$,
\[
E \ll c^M \frac{(cR)^{M+1}}{|w|^{M+1}} \sum_{\mu=0}^M \frac{(cR)^\mu}{\mu!} \ll c^M\frac{(c^2R)^M}{|w|^{\delta M}} \sum_{\mu=0}^M \frac{R^\mu}{\mu!V^{(1-\delta)\mu}} \ll c^M\frac{R^M}{|w|^{\delta M}} e^{O(RV^{\delta-1})},
\]
and the lemma follows. \fine

\medskip
Since $Q_0(s)=1$ identically, writing
\[
\begin{split}
h_M^{**}(w,s) &= \frac{\Gamma(a+\frac12 -\frac{w}{2})}{\Gamma(b+\frac12+\frac{w}{2})} \Gamma(w) \big(1 + \sum_{\nu=1}^M \frac{Q_\nu(s)}{(b+\frac12+\frac{w}{2})\cdots (b+\frac12+\frac{w}{2}+\nu-1)}\big) \\
&= \sum_{\nu=0}^M Q_\nu(s) \frac{\Gamma(a+\frac12 -\frac{w}{2})}{\Gamma(b+\frac12+\frac{w}{2}+\nu)} \Gamma(w) \\
I^{**}_{X.M}(s) &= \frac{1}{2\pi i} \int_{\LL_\infty} h_M^{**}(w,s) (2\omega_X)^{-w} \d w
\end{split}
\]
from \eqref{3-11}, \eqref{3-12} and Lemma 3.6 we obtain
\begin{equation}
\label{3-21}
I_X(s) = c_0(F) \big(\frac{d}{2\beta^{1/d}}\big)^{ds} I^{**}_{X.M}(s) + O(A^{|\si|} e^{O(R^2/V + RV^{\delta-1})} (cR)^{M+2})
\end{equation}
uniformly for $X\geq X_0$, provided $V\geq (c_1R)^{1/\delta}$ and $\max(\frac{2\delta-1}{1-\delta}, \frac{1}{\delta}(d|\si|+d/2+2)) \leq M \leq DR$. Indeed, by Lemma 3.6 we have
\[
h_M^{*}(w,s) = h_M^{**}(w,s) + O\big((cR)^M e^{O(RV^{\delta-1})} \big|\frac{\Gamma(a+\frac12-\frac{w}{2})}{\Gamma(b+\frac12+\frac{w}{2})}\Gamma(w)\big| \frac{1}{|w|^{\delta M}}\big)
\]
and, for $w\in\LL_\infty$, $\Gamma(w) \ll e^{-\frac{\pi}{2}v} |w|^{-K-1/2}$ and $|(2\omega_X)^{-w}| \ll e^{\frac{\pi}{2}v}$. Hence, thanks to \eqref{3-10}, replacing $h_M^{*}(w,s)$ by $h_M^{**}(w,s)$ in $I^*_{X,M}(s)$ causes an error of size
\[
\ll A^{|\si|} (cR)^M e^{O(R^2/V + RV^{\delta-1})} \int_V^\infty v^{-(\delta M-d|\si|-d/2)} \d v \ll A^{|\si|} (cR)^M e^{O(R^2/V + RV^{\delta-1})}
\]
since the integral is convergent under the above conditions on $M$. Therefore, \eqref{3-21} follows from \eqref{3-12}.

\medskip
{\bf Lemma 3.7.}  {\sl Let $Q_\nu(s)$ be as in \eqref{3-20} and $1\leq \nu \leq DR$, where $D\geq 1$ is arbitrary. Then
\[
Q_\nu(s) \leq \frac{(c_7R)^{2\nu}}{\nu!},
\]
where $c_7$ may depend also on $D$.}

\medskip
{\it Proof.} From \eqref{3-17}, \eqref{3-20} and the definition $b=b(s) = \frac{ds}{2} + \frac{\xi_F}{2}$ we have
\[
Q_\nu(s) \ll c^\nu \sum_{\mu=0}^\nu \frac{1}{2^\mu} |V_\mu(s)| \sum_{k=0}^{\nu-\mu} |C_{\mu+k,\nu}| \left|{-\mu\choose k}\right| (R+1)^k.
\]
Recalling that in the proof of Lemma 3.4 we have $V_{\mu,M}(s) = V_\mu(s)$ for $\mu\leq M$ and that $M\leq DR$, from \eqref{3-14} we have
\[
V_\mu(s) \ll \frac{(cR)^{2\mu}}{\mu!}.
\]
Hence by \eqref{3-16} and $|{-\mu\choose k}|={\mu+k-1 \choose k}$ we get
\[
\begin{split}
Q_\nu(s) &\ll c^\nu \sum_{\mu=0}^\nu \sum_{k=0}^{\nu-\mu} \frac{(\nu-1)!}{\mu!(\mu+k-1)!} {\nu-1\choose \mu+k-1} {\mu+k-1 \choose k} (cR)^{2\mu+k} \\ 
&\ll c^\nu \sum_{\mu=0}^\nu \sum_{k=0}^{\nu-\mu} \frac{\nu!}{\mu!(\mu+k)!} 2^\nu 2^{\mu+k-1} (cR)^{2\mu+k} \ll c^\nu \nu!  \sum_{\mu=0}^\nu \sum_{k=0}^{\nu-\mu} \frac{(cR)^{2\mu+k}}{\mu!(\mu+k)!} \\
&\ll c^\nu \nu! \sum_{k=0}^\nu \frac{(cR)^k}{k!} \sum_{\mu=0}^{\nu-k} \left(\frac{(cR)^\mu}{\mu!}\right)^2 \ll c^\nu \nu! \sum_{k=0}^\nu \frac{(cR)^k}{k!} \left(\sum_{\mu=0}^{\nu-k} \frac{(cR)^\mu}{\mu!}\right)^2
\end{split}
\]
thanks to $(\mu+k)!\geq \mu!k!$. Since $\nu\leq DR$ and since the last constant $c$ may clearly be chosen larger than $D$, we may apply Lemma 3.6 of \cite{Ka-Pe/twist}, thus getting (recall that the value of $c$ is not necessarily the same at each occurrence !)
\[
Q_\nu(s) \ll c^\nu \nu! \sum_{k=0}^\nu \frac{(cR)^k}{k!} \big(1+ \frac{(cR)^{2(\nu-k)}}{((\nu-k)!)^2}\big),
\]
and by repeated applications of Lemma 3.6 of \cite{Ka-Pe/twist} we have
\[
\begin{split}
Q_\nu(s) &\ll c^\nu\nu! \left(1 + \frac{(cR)^\nu}{\nu!} + \sum_{k=0}^\nu \frac{(cR)^{2\nu-k}}{k! (\nu-k)!(\nu-k)!}\right) \\
&\ll c^\nu\nu! + (cR)^\nu + (cR)^\nu\sum_{k=0}^\nu {\nu\choose k} \frac{(cR)^{\nu-k}}{(\nu-k)!} \\
&\ll c^\nu\nu! + (cR)^\nu + (cR)^\nu\sum_{k=0}^\nu \frac{(cR)^{\nu-k}}{(\nu-k)!} \\
&\ll c^\nu\nu! + (cR)^\nu + (cR)^\nu \big(1 + \frac{(cR)^\nu}{\nu!}\big) \ll \frac{(cR)^{2\nu}}{\nu!}
\end{split}
\]
since $\nu\leq DR$. Lemma 3.7 is therefore proved. \fine

\medskip
For $z\in\CC$ and integer $\nu\geq 1$ we define
\begin{equation}
\label{3-22}
P_\nu(z) = \prod_{\substack{0\leq j\leq \nu-1 \\ |z+j|\geq 1/2}} (z+j).
\end{equation}

\medskip
{\bf Lemma 3.8.} {\sl With the above notation we have $|P_\nu(z)| \geq 2^{-\nu}(\nu-1)!$}.

\medskip
{\it Proof.} By induction. For $\nu=1$ the right hand side equals $1/2$, while the left hand side is at least $1/2$. Since $\nu= |z-(z+\nu)| \leq |z| + |z+\nu|$, we have $|z|\geq \nu/2$ or $|z+\nu|\geq \nu/2$. Hence
\[
P_{\nu+1}(z) = 
\begin{cases} 
zP_\nu(z+1) & \text{if} \ |z|\geq \nu/2 \\
(z+\nu)P_\nu(z) & \text{if} \ |z+\nu|\geq \nu/2.
\end{cases}
\]
Therefore, for the appropriate $\epsilon=0$ or $1$, from the inductive hypothesis we get
\[
|P_{\nu+1}(z)| \geq \frac{\nu}{2} |P_\nu(z+\epsilon)| \geq \frac{\nu}{2} \frac{(\nu-1)!}{2^\nu} = \frac{\nu!}{2^{\nu+1}},
\]
and the lemma follows. \fine

\medskip
The proof of Theorem 2 is based on the following lemma, which summarizes the results obtained so far in this section. Let (recall the definition of $a$ and $b$ after \eqref{3-4})
\begin{equation}
\label{3-23}
I_{X,\nu}(s) = \frac{1}{2\pi i} \int_{(-K)} \frac{\Gamma(a+\frac12 -\frac{w}{2})}{\Gamma(b+\frac12+\frac{w}{2}+\nu)} \Gamma(w)(2\omega_X)^{-w} \d w.
\end{equation}

\medskip
{\bf Lemma E.} {\sl Let $-L\leq \sigma\leq 3/2$ with $L\geq 3/2$, $|s|\leq R$ with $R\geq 1$, $0<\delta<1$, $\eta\delta>1$, $\delta'=\max(0,2-1/\delta)$, $M=[\frac{d}{\delta}L]+H$ where $H=H(F,\delta,\eta)>0$ is a sufficiently large integer. Then, under the hypotheses of Theorem $2$ and with the notation in \eqref{2-2}, \eqref{3-20} and \eqref{3-23}, we have
\[
F_X(s,\alpha) = c(F) \frac{\overline{a(n_\alpha)}}{n_\alpha^{1-s}} q_F^{-s} (\pi d)^{ds}  \sum_{\nu=0}^M Q_\nu(s) I_{X,\nu}(s) + O\big( A^{L} R^{\eta d L+B} e^{C(R^{\delta'} +R^{(3-\eta\delta)/2})}\big)
\]
uniformly for $X\geq X_0$, where $c(F)\neq 0$.}

\medskip
{\it Proof.} We have to estimate the contribution to $I_X(s)$ of the part $(-K-i\infty, -K+iV)$ of the line of integration in the integrals $I_{X,\nu}(s)$, see \eqref{3-23}. Recalling \eqref{3-22}, by the factorial formula we have
\begin{equation}
\label{3-24}
\Gamma(b+\frac12+\frac{w}{2}+\nu) = \Gamma(b+\frac12+\frac{w}{2}) P_\nu(b+\frac12+\frac{w}{2}) (b+\frac12+\frac{w}{2}+j_0),
\end{equation}
where $1\leq j_0\leq \nu-1$ is such that $|b+\frac12+\frac{w}{2}+j_0|<1/2$, if it exists. If such a $j_0$ does not exist, then the factor $b+\frac12+\frac{w}{2}+j_0$ is not present in \eqref{3-24}. Let $w=-K+iv$ with $v<V$. By Lemma 3.8, the reflection formula for the $\Gamma$ function, the inequality $|\frac{z}{\sin (z-\pi j_0)}| \gg \exp(-|z|)$ for $|z|<1/2$ and the definition of $b$ after \eqref{3-4} we get
\begin{equation}
\label{3-25}
\begin{split}
|\Gamma(b+\frac12+\frac{w}{2}+\nu)| &\gg \frac{1}{|\Gamma(\frac12-b-\frac{w}{2})|} \frac{(\nu-1)!}{2^\nu} \frac{|b+\frac12+\frac{w}{2}+j_0|}{|\sin(\pi(b+\frac12+\frac{w}{2}))|} \\
&\gg \frac{1}{|\Gamma(\frac12-b-\frac{w}{2})|} \frac{\nu!}{\nu2^\nu} e^{-\frac{\pi}{2}|dt+v|}.
\end{split}
\end{equation}
Moreover, recalling also the definition of $a$ after \eqref{3-4}, from Stirling's formula we have
\[
\begin{split}
\Gamma(w) &\ll  e^{-\frac{\pi}{2}|v|} |w|^{-K-1/2} \\
\Gamma(a+\frac12 -\frac{w}{2}) &\ll (1+\frac{|dt+v|}{2})^{d|\si|/2 + B} e^{-\frac{\pi}{4}|dt+v|} \\
\Gamma(\frac12-b-\frac{w}{2}) &\ll (1+\frac{|dt+v|}{2})^{d|\si|/2 + B} e^{-\frac{\pi}{4}|dt+v|}.
\end{split}
\]
Hence, since $|(2\omega_X)^{-w}|=|2\omega_X|^K e^{v\arg\omega_X}$, we obtain
\begin{equation}
\label{3-26}
\begin{split}
\int_{-K-i\infty}^{-K+iV} &\frac{\Gamma(a+\frac12 -\frac{w}{2})}{\Gamma(b+\frac12+\frac{w}{2}+\nu)} \Gamma(w)(2\omega_X)^{-w} \d w \\
&\ll \frac{\nu2^\nu}{\nu!} A^{|\si|} \int_{-\infty}^V \frac{(1+|dt+v|)^{d|\si| + B}}{(1+|v|)^{K+1/2}} e^{-\frac{\pi}{2}|v| + v\arg\omega_X} \d v \\
&\ll \frac{\nu2^\nu}{\nu!} A^{|\si|} V^{d|\si| + B} + \frac{\nu2^\nu}{\nu!} A^{|\si|} \int_V^{+\infty} e^{-v} v^{d|\si| + B} \d v \\
&\ll \frac{\nu2^\nu}{\nu!} A^{|\si|} V^{d|\si| + B} \big\{1 + e^{-V}(1+\frac{1+|\si|^2}{V})\big\} \ll \frac{\nu2^\nu}{\nu!} A^{|\si|} V^{d|\si| + B} 
\end{split}
\end{equation}
thanks to \eqref{2-16} and the fact that $V\geq (c_1R)^{1/\delta}$. Thanks to \eqref{3-26}, to the choice of $M$ and to Lemma 3.7, the contribution of the left hand side of \eqref{3-26} to $I_X(s)$ is
\begin{equation}
\label{3-27}
\begin{split}
&\ll A^{|\si|} V^{d|\si| + B} M2^M \sum_{\nu=0}^M \frac{|Q_\nu(s)|}{\nu!} \ll A^{L} V^{d|\si| + B} \sum_{\nu=0}^M \frac{(CR)^{2\nu}}{(\nu!)^2} \\
&\ll A^{L} V^{d|\si| + c} \sum_{\nu=0}^M \frac{(CR)^{(3-\eta\delta)\nu+(\eta\delta-1)\nu}}{(\nu!)^2} \ll A^{L} V^{d|\si| + c} (cR)^{(\eta\delta-1)M} \big(\sum_{\nu=0}^\infty \frac{(CR)^{\frac{3-\eta\delta}{2}\nu}}{\nu!}\big)^2 \\
&\ll A^{L} V^{d|\si| + c} (cR)^{(\eta\delta-1)M} e^{CR^{\frac{3-\eta\delta}{2}}}.
\end{split}
\end{equation}
Recalling that $q_F=(2\pi)^dQ^2\beta$ and writing $c(F) = \omega Q c_0(F)$, Lemma E follows now from \eqref{3-1},\eqref{3-2},\eqref{3-21},\eqref{3-23} and \eqref{3-27}, choosing $V=(cR)^{1/\delta}$. \fine

\bigskip
Now we show how Theorem 2 follows from Lemma E and some of the arguments in our previous treatment of the standard twist in \cite{Ka-Pe/1999a} and \cite{Ka-Pe/2005}. Under the hypotheses of Lemma E and with the notation there we write
\[
H_X(s,\alpha) = F_X(s,\alpha) - c(F) \frac{\overline{a(n_\alpha)}}{n_\alpha^{1-s}} q_F^{-s} (\pi d)^{ds}  \sum_{\nu=0}^M Q_\nu(s) I_{X,\nu}(s)
\]
and note that for $\si>1$
\[
\lim_{X\to\infty} H_X(s,\alpha) = H(s,\alpha) 
\]
exists. Indeed, clearly $F_X(s,\alpha)$ tends to $F(s,\alpha)$ for $\si>1$, while the treatment of the limit of the integrals $I_{X,\nu}(s)$ as $X\to\infty$ is borrowed from our previous papers on this subject, see in particular Theorem 5.1 of \cite{Ka-Pe/1999a} in the case $d=1$ and Lemma 2.4 of \cite{Ka-Pe/2005} in the general case. Note, however, that Lemma 2.4 of \cite{Ka-Pe/2005} deals with a normalized situation where $s$, $\lambda_j$ and $\mu_j$ mean, in our present notation,
\[
ds - \frac{d-1}{2}, \hskip1cm \frac{\lambda_j}{d}, \hskip1cm \mu_j + \frac{\lambda_j}{2}(1-\frac{1}{d}),
\]
respectively (see the discussion before and after $(2.2)$ and $(2.3)$ on pages 320-321 of \cite{Ka-Pe/2005}). As a consequence we have
\begin{equation}
\label{3-28}
H(s,\alpha) = F(s,\alpha) - c(F) \frac{\overline{a(n_\alpha)}}{n_\alpha^{1-s}} q_F^{-s} (\pi d)^{ds}  \sum_{\nu=0}^M Q_\nu(s) I_{\nu}(s),
\end{equation}
where, recalling that $a=a(s) = \frac{d}{2}(1-s)+\frac12\overline{\xi_F}$ and $b=b(s) = \frac{d}{2}s+\frac12\xi_F$ (see after \eqref{3-4}),
\begin{equation}
\label{3-29}
\begin{split}
I_{\nu}(s) &= -\sum_{k=0}^K \frac{1}{k!} \frac{\Gamma(a+\frac12+\frac{k}{2})}{\Gamma(b+\frac12-\frac{k}{2}+\nu)}
(\frac{i}{2})^{-k} +\sqrt{\pi} \frac{\Gamma(a+\frac12)\Gamma(b-a-\frac12+\nu)} {\Gamma(-a)\Gamma(b+\nu)\Gamma(\frac12+b+\nu)} \\
&-i\sqrt{\pi} \frac{\Gamma(a+1)\Gamma(b-a-\frac12+\nu)} {\Gamma(\frac12-a)\Gamma(b+\nu)\Gamma(\frac12+b+\nu)} = -I_\nu^{(1)}(s) + \sqrt{\pi}I_\nu^{(2)}(s) -i\sqrt{\pi}I_\nu^{(3)}(s),
\end{split}
\end{equation}
say. Theorem 5.1 of \cite{Ka-Pe/1999a} deals only with $I_0(s)$, while Lemma 2.4 of \cite{Ka-Pe/2005} deals with all $I_\nu(s)$, there denoted by $\Gamma_{K,\nu}(s)$. From Lemma E and Vitali's convergence theorem, see Lemma C and the proof of Theorem 1 in Section 2, we therefore have that the limit function $H(s,\alpha)$ exists and is holomorphic for $|s|\leq R$ and $-L\leq \si\leq 3/2$, and satisfies
\begin{equation}
\label{3-30}
H(s,\alpha) \ll A^{L} R^{\eta d L+B} e^{C(R^{\delta'} +R^{(3-\eta\delta)/2})}.
\end{equation}
Since $R$ and $L$ are arbitrary, from \eqref{3-28} and \eqref{3-29} we have in particular that $F(s,\alpha)$ has meromorphic continuation to $\CC$. Moreover, from Lemma 2.5 of \cite{Ka-Pe/2005} (recall the normalization there) we have that $F(s,\alpha)$ is holomorphic over $\CC$ apart possibly from simple poles at the points $s_k$ in \eqref{1-2} satisfying \eqref{1-3}, coming from the term $\Gamma(b-a-\frac12+\nu)$ in \eqref{3-29}. Our next goal is therefore estimating the last term in \eqref{3-28} away from such poles.

\medskip
Before starting the estimation of such a term, we remark that we shall always assume that $s$ is $\delta_0$-apart from any pole which might arise during the estimation. We call potential poles such poles, and we shall deal with them at the end of the proof. Moreover, we always assume that $0\leq k\leq K$ and $|s|\leq R$. We deal first with $I_\nu^{(1)}(s)$, in a similar way as in the proof of Lemma E. Indeed, following \eqref{3-24} we have
\[
\Gamma(b+\frac12-\frac{k}{2}+\nu) = \Gamma(b+\frac12-\frac{k}{2}) P_\nu(b+\frac12-\frac{k}{2}) (b+\frac12-\frac{k}{2}-j_0(k)),
\]
where $1\leq j_0(k)\leq \nu-1$ is such that $|b+\frac12-\frac{w}{2}-j_0(k)|<1/2$, if it exists. Hence similarly as in \eqref{3-25} we get
\[
|\Gamma(b+\frac12-\frac{k}{2}+\nu)| \gg \frac{1}{|\Gamma(\frac12-b+\frac{k}{2})|} \frac{\nu!}{\nu2^\nu} e^{-\frac{\pi}{2}d|t|}
\]
and hence
\[
I_\nu^{(1)}(s) \ll \frac{\nu2^\nu}{\nu!} \sum_{k=0}^K \frac{1}{k!} |\Gamma(a+\frac12+\frac{k}{2}) \Gamma((\frac12 -b +\frac{k}{2})| e^{\frac{\pi}{2}d|t|}.
\]
Using \eqref{2-5}, away from the potential poles we have
\[
|\Gamma(a+\frac12+\frac{k}{2}) \Gamma(\frac12 - b +\frac{k}{2})| \ll A^{L} R^{dL+B} e^{-\frac{\pi}{2}d|t|},
\]
therefore under the above mentioned assumptions we have
\begin{equation}
\label{3-31}
I_\nu^{(1)}(s) \ll \frac{A^{L}}{\nu!} R^{dL+B} .
\end{equation}
In order to estimate $I_\nu^{(2)}(s)$ we apply the duplication formula, the factorial formula and recall \eqref{3-22}, thus getting
\[
\begin{split}
\Gamma(b+\nu)&\Gamma(\frac12 + b + \nu) = \frac{2^{-ds-\xi_F-2\nu}}{\sqrt{\pi}} \Gamma(ds+\xi_F+2\nu) \\
&=\frac{2^{-ds-\xi_F-2\nu}}{\sqrt{\pi}} \Gamma(ds+\xi_F+\nu) P_\nu(ds+\xi_F+\nu) (ds+\xi_F+j_0(\nu))
\end{split}
\]
where $\nu\leq j_0(\nu)\leq 2\nu-1$ is such that $|ds+\xi_F+j_0(\nu)|<1/2$, if it exists. Hence by the reflection formula we obtain
\begin{equation}
\label{3-32}
\begin{split}
I_\nu^{(2)}(s) = &\sqrt{\pi}2^{ds+\xi_F+2\nu} \frac{\sin(\pi(\frac{ds}{2}-\frac{d}{2}-\frac{\overline{\xi_F}}{2}))}{ (ds+\xi_F+j_0(\nu)) P_\nu(ds+\xi_F+\nu)} \times \\
&\times \frac{\Gamma(\frac{d+1}{2}-\frac{ds}{2}+\frac{\overline{\xi_F}}{2}) \Gamma(1+\frac{d}{2}-\frac{ds}{2}+\frac{\overline{\xi_F}}{2}) \Gamma(ds-\frac{d+1}{2}+id\theta_F+\nu)}{\Gamma(ds+\xi_F+\nu)},
\end{split}
\end{equation}
and by Lemma 3.8 we get, away from potential poles,
\[
\begin{split}
I_\nu^{(2)}(s) \ll &\frac{A^{L} c^\nu}{\nu!} e^{\frac{\pi}{2}d|t|} |\Gamma(\frac{d+1}{2}-\frac{ds}{2}+\frac{\overline{\xi_F}}{2}) \Gamma(1+\frac{d}{2}-\frac{ds}{2}+\frac{\overline{\xi_F}}{2})| \times \\
&\times \left| \frac{\Gamma(ds-\frac{d+1}{2}+id\theta_F)}{\Gamma(ds+\xi_F)}\right| \prod_{j=0}^{\nu-1} \left| \frac{ds-\frac{d+1}{2}+id\theta_F+j}{ds+\xi_F+j}\right|.
\end{split}
\]
Since each factor of the last product is bounded away from the potential poles we have
\[
\prod_{j=0}^{\nu-1} \left| \frac{ds-\frac{d+1}{2}+id\theta_F+j}{ds+\xi_F+j}\right| \ll c^\nu,
\]
and by Stirling's formula 
\[
\left| \Gamma(\frac{d+1}{2}-\frac{ds}{2}+\frac{\overline{\xi_F}}{2}) \Gamma(1+\frac{d}{2}-\frac{ds}{2}+\frac{\overline{\xi_F}}{2})\right| \ll e^{-\frac{\pi}{2}d|t|} R^{dL+B}.
\]
Moreover, by the reflection formula and Stirling's formula we get
\[
\left| \frac{\Gamma(ds-\frac{d+1}{2}+id\theta_F)}{\Gamma(ds+\xi_F)}\right| \ll \left|\frac{\Gamma(1-ds-\xi_F)}{\Gamma(\frac{d+3}{2}-ds-id\theta_F)} \right| \ll R^B.
\]
Therefore, from \eqref{3-32} and the above bounds we conclude that under the above assumptions
\begin{equation}
\label{3-33}
I_\nu^{(2)}(s) \ll \frac{c^\nu}{\nu!} A^{L} R^{dL+B}.
\end{equation}
Similar computations show that the same bound holds for $I_\nu^{(3)}(s)$ as well, hence from \eqref{3-29}, \eqref{3-31} and \eqref{3-33} we obtain
\[
I_\nu(s) \ll \frac{c^\nu}{\nu!} A^{L} R^{dL+B}
\]
for $0\leq k\leq K$ and $|s|\leq R$ away from all the above potential poles. Therefore from Lemma 3.7 we have
\[
\begin{split}
c(F) \frac{\overline{a(n_\alpha)}}{n_\alpha^{1-s}} q_F^{-s} (\pi d)^{ds}  \sum_{\nu=0}^M Q_\nu(s) I_{\nu}(s) &\ll  A^{L} R^{dL+B} \sum_{\nu=0}^M \frac{(CR)^{2\nu}}{(\nu!)^2} \\
&\ll A^{L} R^{\eta dL+B} \sum_{\nu=0}^M \frac{(CR)^{2\nu-(\eta-1)dL}}{(\nu!)^2},
\end{split}
\]
and since $(\eta-1)dL \geq (\eta-1)(\delta M-c) \geq \delta(\eta-1)\nu - c$ we obtain
\[
\begin{split}
\sum_{\nu=0}^M \frac{(CR)^{2\nu-(\eta-1)dL}}{(\nu!)^2} &\ll R^c \sum_{\nu=0}^M \frac{(cR)^{(2-(\eta-1)\delta)\nu}}{(\nu!)^2} \ll R^c \left(\sum_{\nu=0}^\infty \frac{(CR)^{(1-(\eta-1)\delta/2)\nu}}{(\nu!)^2}\right)^2\\ 
&\ll R^ce^{CR^{1-\frac{(\eta-1)\delta}{2}}} \ll R^c e^{CR^{(3-\eta\delta)/2}},
\end{split}
\]
in view of $1-(\eta-1)\delta/2 \leq (3-\eta\delta)/2$. As a consequence, away from all the above potential poles we have
\begin{equation}
\label{3-34}
c(F) \frac{\overline{a(n_\alpha)}}{n_\alpha^{1-s}} q_F^{-s} (\pi d)^{ds}  \sum_{\nu=0}^M Q_\nu(s) I_{\nu}(s) \ll A^{L} R^{\eta dL+B} e^{CR^{(3-\eta\delta)/2}}.
\end{equation}
Hence from \eqref{3-28}, \eqref{3-30} and \eqref{3-34} we finally obtain
\begin{equation}
\label{3-35}
F(s,\alpha) \ll A^{L} R^{\eta d L+B} e^{C(R^{\delta'} +R^{(3-\eta\delta)/2})}
\end{equation}
away from the potential poles, provided $|s|\leq R$ and $-L\leq \sigma\leq 3/2$. However, we already know from Theorem 1 of \cite{Ka-Pe/2005} that $F(s,\alpha)$ is holomorphic over $\CC$ apart possibly at $s=s_k$, $k=0,1,\dots$. Hence, thanks to the maximum modulus principle, the bound \eqref{3-35} holds for every $s$ outside discs of radius $\delta_0$ around the points $s_k$, and Theorem 2 follows.

\bigskip
\section{Proof of Theorem 4}

\smallskip
Let $f(n,\bfalpha)$ be as in Theorem 4 and define
\[
\A_f = \{\omega=\sum_{\nu=0}^Nm_\nu\kappa_\nu: m_\nu\in\ZZ, m_\nu\geq 0\}
\]
and
\[
A(\omega) = \sum_{\substack{m_0\geq 0,\dots,m_N\geq0 \\ \sum_{\nu=0}^Nm_\nu\kappa_\nu = \omega}} \frac{\prod_{\nu=0}^N(-2\pi i\alpha_\nu)^{m_\nu}}{m_0!\cdots m_N!};
\]
note that $\A_f$ has no accumulation points in $\RR$. We have

\medskip
{\bf Lemma 4.1.} {\sl We have the bounds
\[
\sum_{\substack{\omega \in\A_f \\ \omega\leq x}}1 \ll x^{N+1} \quad \text{and} \quad A(\omega) \ll e^{-c\omega\log\omega} \quad \text{with} \quad c=c(f)>0,
\]
and the series $\sum_{\omega\in\A_f} A(\omega)$ is absolutely convergent.}

\medskip
{\it Proof.} The first bound is trivial, since for every $\nu$ we have $m_\nu\leq 1 + x/\kappa_\nu$ and hence the number of possible $\omega\in\A_f$ up to $x$ is $\ll \prod_{\nu=0}^N (1+x/\kappa_\nu) \ll x^{N+1}$. To prove the second bound we write $r(\omega)$ for the number of representations $\sum_{\nu=0}^Nm_\nu\kappa_\nu = \omega$ and $M_\nu$ for the maximum of the $m_\nu$ occurring in such representations. Recalling that the value of the constants $c$ may not be the same at each occurrence we have
\[
A(\omega) \ll \frac{r(\omega) \prod_{\nu=0}^N (2\pi|\alpha_\nu|+1)^{M_\nu}}{\big[\frac{\omega}{\kappa_N}\big]!}
\ll \frac{\omega^{N+1} \prod_{\nu=0}^N (2\pi|\alpha_\nu|+1)^{\omega/\kappa_\nu}}{\big[\frac{\omega}{\kappa_N}\big]!} \ll \frac{\omega^{N+1} c^\omega}{\big[\frac{\omega}{\kappa_N}\big]!} \ll \frac{c^\omega}{\omega^{\omega/\kappa_N}},
\]
and the second bound follows. Finally, the absolute convergence of the series follows by partial summation from the first two bounds. \fine

\medskip
Let now $F(s)$ be as in Theorem 4. Then, expanding the exponential, thanks to the good convergence properties of the involved series we get for $\si>1$
\begin{equation}
\label{T4-1}
\begin{split}
F(s;f) &= \sum_{n=1}^\infty \frac{a(n)}{n^s} e(-f(n,\bfalpha)) =  \sum_{n=1}^\infty \frac{a(n)}{n^s} \prod_{\nu=0}^N e^{-2\pi i\alpha_\nu n^{-\kappa_\nu}} \\
&=  \sum_{n=1}^\infty \frac{a(n)}{n^s} \sum_{m_0\geq0,\dots,m_N\geq0} \frac{\prod_{\nu=0}^N (-2\pi i\alpha_\nu n^{-\kappa_\nu})^{m_\nu}}{m_0!\cdots m_N!} \\
&= \sum_{m_0\geq0,\dots,m_N\geq0} \frac{\prod_{\nu=0}^N (-2\pi i\alpha_\nu)^{m_\nu}}{m_0!\cdots m_N!} \sum_{n=1}^\infty \frac{a(n)}{n^{s+\sum_{\nu=0}^N m_\nu\kappa_\nu}} \\
&= \sum_{\omega\in\A_f} A(\omega)F(s+\omega),
\end{split}
\end{equation}
the last series being clearly absolutely convergent for $\si>1$ thanks to Lemma 4.1. Moreover, for any given $s\in\CC$ we have $\si+\omega\geq 2$ apart from finitely many values of $\omega\in\A_f$, hence such a series is absolutely and uniformly convergent on compact subsets of $\CC$, apart from the points $s=s_0-\omega$ where $s_0$ is a pole of $F(s)$ and $\omega\in\A_f$ is such that $A(\omega)\neq 0$. Therefore $F(s;f)$ is meromorphic on $\CC$, its poles are contained in the same horizontal strip of $F(s)$, and $F(s;f)$ is entire if $F(s)$ is entire. Further, let $\si$ be fixed and note that by the above observation we may clearly write
\begin{equation}
\label{T4-2}
F(s;f) = F(s) + \sum_{\substack{\omega\in\A_f \\ 0<\omega\leq \omega_0}} A(\omega)F(s+\omega) + B(s)
\end{equation}
where $\omega_0=\omega_0(\si)$ is such that $\si+\omega_0\leq 3/2$ and $B(s)$ is bounded as $|t|\to\infty$. Thanks to the properties of the Lindel\"of $\mu$-function of Dirichlet series, in particular the fact that it is strictly decreasing on any interval where it is positive, \eqref{T4-2} immediately shows that $\mu_F(\si)=\mu_F(\si;f)$, and Theorem 4 is proved.

\bigskip
\section{Proof of Theorem 5}

\smallskip
Suppose first that $F\in M(\rho,\tau)$, $\tau>0$, and $0<\lambda<1/\rho$. As in Sections 2 and 3 we may assume that $\alpha>0$ and $\si\leq 3/2$; we deal first with the case $t>0$. Let $u_0=(2-\si)/\lambda$, $v_0=\tau/\lambda$,
\[
\LL_{-\infty} = (u_0-i\infty, u_0+iv_0), \hskip1.5cm \LL_{v_0}=(u_0+iv_0,-\infty+iv_0)
\]
and, using the notation in the proof of Theorem 1, write
\[
F^\lambda_X(s,\alpha) = \sum_{n=1}^\infty\frac{a(n)}{n^s}e^{-z_Xn^\lambda}.
\]
Since $F(s)$ is holomorphic for $|t|\geq \tau$, and $\Im(s+\lambda w)>\tau$ when $v\geq v_0$, by Mellin's transform and arguments already used in the proof of Theorem 1 we have
\begin{equation}
\label{T5-1}
\begin{split}
F^\lambda_X(s,\alpha) &= \frac{1}{2\pi i}\int_{(u_0)} F(s+\lambda w) \Gamma(w) z_X^{-w}\d w \\
&=  \frac{1}{2\pi i}\left(\int_{\LL_{-\infty}} + \int_{\LL_{v_0}}\right) F(s+\lambda w) \Gamma(w) z_X^{-w}\d w = F^{(1)}_X(s) + F^{(2)}_X(s),
\end{split}
\end{equation}
say. Clearly,
\[
F^{(1)}_X(s)  \ll \int_{-\infty}^{v_0} |\Gamma(u_0+iv)z_X^{-u_0-iv}|\d v
\]
and for $w\in\LL_{-\infty}$ we have
\[
\begin{split}
\log|\Gamma(w)| &= \Re\big((w-\frac12)\log w-w\big)+ O(1) = (u_0-\frac12)\log |w| -v\arctan\frac{\lambda v}{2-\si} + O(1+|\si|) \\
&= (u_0-\frac12)\log |w| -|v|\big(\frac{\pi}{2} + O(\frac{|\si|+1}{|v|})\big) +O(1+|\si|)),
\end{split}
\]
therefore
\[
\Gamma(w) \ll A^{|\si|} |w|^{u_0-1/2} e^{-\frac{\pi}{2}|v|}.
\]
Moreover
\[
|z_X^{-w}| = |z_X|^{-u_0} e^{v\arg z_X} \ll A^{|\si|}e^{v(\frac{\pi}{2}+O(1/X))},
\]
hence
\begin{equation}
\label{T5-2}
\begin{split}
F^{(1)}_X(s) &\ll A^{|\si|} \int_{-\infty}^{v_0} |w|^{u_0-1/2} e^{-\frac{\pi}{2}|v| + \frac{\pi}{2}v+O(\frac{|v|}{X})} \d v \\
&\ll A^{|\si|} + A^{|\si|} \int_0^\infty |u_0+iv|^{u_0-1/2} e^{-v(\pi+O(1/X))} \d v \\
&\ll A^{|\si|} + A^{|\si|}(1+|\si|)^{|\si|/\lambda + B} + A^{|\si|} \int_0^\infty (1+v)^{u_0-1/2}  e^{-v(\pi+O(1/X))} \d v \\
&\ll A^{|\si|} + A^{|\si|}(1+|\si|)^{|\si|/\lambda + B} + A^{|\si|} \Gamma(\frac{|\si|}{\lambda} +B) \ll A^{|\si|} (1+|s|)^{|\si|/\lambda + B}
\end{split}
\end{equation}
uniformly for $X\geq X_0$.

\medskip
Recalling \eqref{1-6}, for $w\in\LL_{v_0}$ we have
\[
\begin{split}
F(s+\lambda w) &\ll A^{|\si+\lambda u|} (1+|s+\lambda w|)^{\rho|\si+\lambda u|+B} e^{C|s+\lambda w|^\delta} \\
&\ll A^{|\si| +|\lambda u|} (1+|s+\lambda w|)^{\rho|\si| - \lambda \rho u +B} e^{C(|s|^\delta+|u|^\delta)},
\end{split}
\]
and also, thanks to the reflection formula and Stirling's formula,
\[
\Gamma(w) \ll \frac{1}{|\Gamma(1-w)|} \ll A^{|u|} (1+|u|)^{-|u|}, \hskip2cm z_X^{-w} \ll A^{|u|}.
\]
Therefore
\begin{equation}
\label{T5-3}
\begin{split}
F^{(2)}_X(s) &\ll A^{|\si|} e^{C|s|^\delta} \int_{-\infty}^{u_0} A^{|u|} (1+|s+\lambda w|)^{\rho|\si| - \lambda \rho u +B} (1+|u|)^{-|u|} \d u \\
&\ll A^{|\si|} e^{C|s|^\delta} \int_{-\infty}^{-|s|} A^{|u|} (1+|u|)^{-(1-\lambda\rho)|u| + \rho|\si|+B} \d u \\
&\hskip1.5cm + A^{|\si|} e^{C|s|^\delta} (1+|s|)^{\rho|\si|+B} \int_{-|s|}^{u_0} A^{|u|} \frac{(1+|s|)^{-\lambda \rho u}}{(1+|u|)^{|u|}} \d u \\
&\ll A^{|\si|} e^{C|s|^\delta} (1+|s|)^{\rho|\si|+B} \big(1+ \int_1^{|s|} A^{u} \frac{(1+|s|)^{\lambda \rho u}}{u^u} \d u \big) \\
&\ll A^{|\si|} (1+|s|)^{\rho|\si|+B} e^{C(|s|^\delta +|s|^{\lambda \rho})}
\end{split}
\end{equation}
uniformly for $X\geq X_0$, since the integrand in the last but one row is maximal (essentially) when $u=c|s|^{\lambda \rho}$. Hence from \eqref{T5-1}, \eqref{T5-2} and \eqref{T5-3} we obtain that for $\sigma\leq 3/2$ and $t>0$
\begin{equation}
\label{T5-4}
F^\lambda_X(s,\alpha) \ll A^{|\si|} (1+|s|)^{|\si|/\lambda+B} e^{C(|s|^\delta +|s|^{\lambda \rho})}
\end{equation}
uniformly for $X\geq X_0$.

\medskip
The case $t\leq0$ is partly similar to the previous case, thus we will be more synthetic; all the estimates below hold uniformly for $X\geq X_0$. Given $F\in M(\rho,\tau)$ we write $u_0=(2-\si)/\lambda$, $v_0=(\tau-t)/\lambda$ and consider again the half-lines $\LL_{-\infty}$ and $\LL_{v_0}$ (with the new value of $v_0$). Since $\Im(s+\lambda w) \geq\tau$ for $v\geq v_0$, as in the previous case we have
\begin{equation}
\label{T5-5}
\begin{split}
F^\lambda_X(s,\alpha) =  \frac{1}{2\pi i}\left(\int_{\LL_{-\infty}} + \int_{\LL_{v_0}}\right) F(s+\lambda w) \Gamma(w) z_X^{-w}\d w = F^{(3)}_X(s) + F^{(4)}_X(s),
\end{split}
\end{equation}
say, and a similar argument gives, uniformly for $X\geq X_0$,
\begin{equation}
\label{T5-6}
F^{(3)}_X(s) \ll A^{|\si|} (1+|s|)^{|\si|/\lambda + B}.
\end{equation}

\medskip
Recalling again \eqref{1-6}, for $w\in\LL_{v_0}$ we have
\[
F(s+\lambda w) \ll A^{|\si| +|\lambda u|} (1+|\si+\lambda u|)^{\rho|\si + \lambda u| +B} e^{C|\si +\lambda u|^\delta}.
\]
Moreover, for $w\in\LL_{v_0}$ with $u\geq0$ from Stirling's formula we get
\[
\begin{split}
|\Gamma(w)z_X^{-w}| &\ll A^{|\si|} (1+|s|)^u e^{-u} e^{-v_0\arg w} e^{v_0\arg z_X} \\
&\ll A^{|\si|} (1+|s|)^u e^{-v_0(u/v_0 + \arctan v_0/u)} e^{\pi v_0/2} \ll A^{|\si|} (1+|s|)^u
\end{split}
\]
since $x+\arctan 1/x\geq \pi/2$ for $x>0$. Hence the contribution to $F^{(4)}_X(s)$ coming from the part of $\LL_{v_0}$ with $u\geq 0$ is
\begin{equation}
\label{T5-7}
\ll A^{|\si|} \int_0^{u_0} (1+|s|)^{\rho|\si|+u(1-\lambda \rho) +B} \d u \ll A^{|\si|} (1+|s|)^{|\si|/\lambda + B}.
\end{equation}
For $w\in\LL_{v_0}$ with $u<0$ we have, again by the reflection and Stirling formulae, that
\[
|\Gamma(w)z_X^{-w}| \ll e^{-\pi v_0} \frac{1}{|\Gamma(1-w)|} A^{|u|} e^{\pi v_0/2} \ll |w|^u A^{|u|}.
\]
Hence the contribution to $F^{(4)}_X(s)$ coming from the part of $\LL_{v_0}$ with $-c(1+|\si|)\leq u \leq0$ (here $c>1$ is an arbitrary constant) is
\begin{equation}
\label{T5-8}
\begin{split}
&\ll A^{c|\si|} (1+|s|)^{\rho|\si|+B} \int_0^{c(1+|\si|)} \frac{(1+|\si|)^{\lambda \rho u}}{|w|^u} \d u \\
&\ll A^{c|\si|} (1+|s|)^{\rho|\si|+B} \max_{0\leq u\leq c(1+|\si|)} \frac{(1+|\si|)^{\lambda \rho u}}{|w|^u} \ll A^{c|\si|} (1+|s|)^{\rho|\si|+B},
\end{split}
\end{equation}
 since (recall that $\lambda < 1/\rho$)
 \[
 \max_{0\leq u\leq c(1+|\si|)} \left(\frac{(1+|\si|)^{\lambda \rho}}{|w|}\right)^u \leq 1+  \max_{0\leq u\leq (1+|\si|)^{\lambda \rho}} \left(\frac{(1+|\si|)^{\lambda \rho}}{|w|}\right)^u \ll (1+|\si|)^{\lambda \rho(1+|\si|)^{\lambda \rho}} \ll A^{|\si|}.
 \]
 Finally, the contribution to $F^{(4)}_X(s)$ coming from the part of $\LL_{v_0}$ with $u<-c(1+|\si|)$ is
 \[
 \ll A^{|\si|} \int_{c(1+|\si|)}^\infty A^u u^{\rho|\si| + (\lambda \rho-1)u+B} \d u = A^{|\si|} \int_{c(1+|\si|)}^\infty e^{f(u)}\frac{\d u}{u^2},
 \]
 where $f(u) = u\log A + (\rho|\si|+B)\log u - (1-\lambda \rho)u\log u +2$. But $f'(u) = -(1-\lambda \rho)\log u +O(1)$, therefore $f'(u)<0$  for $u>c(1+|\si|)$ if $c$ is large enough and hence such a contribution is
\begin{equation}
\label{T5-9}
\ll A^{|\si|} e^{f(c(1+|\si|))} \ll A^{|\si|} (1+|s|)^{\rho|\si|+B},
\end{equation}
with such a choice of $c$. From \eqref{T5-7}-\eqref{T5-9} we obtain that 
\begin{equation}
\label{T5-10}
F^{(4)}_X(s) \ll A^{|\si|} (1+|s|)^{|\si|/\lambda +B},
\end{equation}
hence from \eqref{T5-5}, \eqref{T5-6} and \eqref{T5-10} we have that for $\si\leq 3/2$ and $t\leq0$
\begin{equation}
\label{T5-11}
F^\lambda_X(s,\alpha) \ll A^{|\si|} (1+|s|)^{|\si|/\lambda +B}
\end{equation}
uniformly for $X\geq X_0$.

\medskip
From \eqref{T5-4} and \eqref{T5-11} we get that for $\si\leq 3/2$
\begin{equation}
\label{T5-12}
F^\lambda_X(s,\alpha) \ll A^{|\si|} (1+|s|)^{|\si|/\lambda+B} e^{C(|s|^\delta +|s|^{\lambda \rho})}
\end{equation}
uniformly for $X\geq X_0$. Hence the argument based on Vitali's convergence theorem in the proof of Theorem 1 shows that $F^\lambda(s,\alpha)$ is entire and satisfies the bound in \eqref{T5-12}, therefore $F^\lambda(s,\alpha)$ belongs to $M(1/\lambda)$ and the first part of Theorem 5 is proved.

\bigskip
Suppose now that $F\in M(\rho)$ and $0<\lambda<1/\rho$; again we may assume that $\alpha>0$. In this case $F(s)$ is entire, and the argument is simpler. Starting with the usual integral representation of $F^{\lambda}_X(s,\alpha)$ (recall that $X$ is a large integer) and shifting the line of integration to $u= -X-1/2$, thanks to the decay of the $\Gamma$ function on vertical strips we obtain
\begin{equation}
\label{T5-13}
\begin{split}
F^\lambda_X(s,\alpha) &= \frac{1}{2\pi i}\int_{(u_0)} F(s+\lambda w) \Gamma(w) z_X^{-w}\d w \\
&= \sum_{k=0}^X \frac{(-1)^k}{k!} F(s-\lambda k) z_X^k + O\big(\int_{(-X-1/2)} |F(s+\lambda w)\Gamma(w) z_X^{-w} \d w|\big).
\end{split}
\end{equation}
Next we show that, as $X\to\infty$, the last integral tends to $0$ and the resulting series is suitably convergent, thus it represents an entire function. This provides the analytic continuation and series representation of $F^\lambda(s,\alpha)$ over the whole complex plane, since $F^\lambda_X(s,\alpha)\to F^\lambda(s,\alpha)$ for $\si>1$. Moreover, we get suitable bounds for such a series, showing that $F^\lambda(s,\alpha)$ belongs to $M(\rho)$ and thus closing the proof of Theorem 5.

\medskip
Once again thanks to \eqref{1-6}, for $u=-X-1/2$ we have 
\[
F(s+\lambda w) \ll A^{|\si|+\lambda(X+1/2)} (|s|+\lambda|w|)^{\rho(|\si|+\lambda X)+B} e^{C(|s|+\lambda|w|)^\delta},
\]
and moreover (reflection and Stirling formulae)
\[
\Gamma(w) \ll e^{-\pi |v|} \frac{1}{|\Gamma(1-w)|} \ll e^{-\pi |v| +|v||\arg(-w)|} |w|^{-X-1} e^X
\]
\[
z_X^{-w} \ll |z_X|^{X+1/2} e^{v\arg z_X}.
\]
Hence
\[
F(s+\lambda w)\Gamma(w)z_X^{-w} \ll A^{|\si|+X} (|s|+\lambda|w|)^{\rho(|\si|+\lambda X)+B} |w|^{-X-1} e^{C(|s|+\lambda|w|)^\delta} e^{-c(X)|v|}
\]
where $c(X)\gg 1/X$. Let $K\subset\CC$ be compact, $s\in K$ and $X\geq X_0(K,F,\rho,\lambda)$. Then the above bound becomes
\[
F(s+\lambda w)\Gamma(w)z_X^{-w} \ll A^X |w|^{-(1-\lambda \rho)X+B} e^{C|w|^\delta -c(X)|v|} \ll A^X X^{-(1-\lambda \rho)X} e^{C|v|^\delta - c|v|/X}
\]
uniformly for $s\in K$ (constants may now depend also on $K$), and therefore
\[
\int_{(-X-1/2)} |F(s+\lambda w)\Gamma(w) z_X^{-w} \d w| \ll A^X X^{-(1-\lambda \rho)X} \to 0
\]
as $X\to\infty$ since $\lambda<1/\rho$. Finally, we show that the series
\begin{equation}
\label{new1}
 \sum_{k=0}^\infty \frac{(-1)^k}{k!} F(s-\lambda k) (2\pi i\alpha)^k,
\end{equation}
obtained from \eqref{T5-13} as $X\to\infty$, is absolutely and uniformly convergent for $s\in K$ and, at the same time, we give the required bounds to show that $F^\lambda(s,\alpha)$ belongs to $M(\rho)$. Indeed, from \eqref{1-6} we get
\begin{equation}
\label{T5-14}
\sum_{k=0}^\infty \frac{1}{k!} |F(s-\lambda k)| (2\pi\alpha)^k \ll \sum_{k=0}^\infty \frac{1}{k!} A^{|\si|+\lambda k} (|s|+\lambda k)^{\rho(|\si|+\lambda k)+B} e^{C(|s|+\lambda k)^\delta} (2\pi\alpha)^k, 
\end{equation}
which shows the absolute and uniform convergence since the term $k!$ dominates over the terms in the numerator, again thanks to the fact that $\lambda<1/\rho$. Moreover, for $|s|$ sufficiently large we split the series in the right hand side of \eqref{T5-14} into $S_1 + S_2$, where $S_1$ is the sum with $k\leq |s|/\lambda$ and $S_2$ with $k> |s|/\lambda$. But since $\lambda<1/\rho$
\[
\begin{split}
S_1 &\ll A^{|\si|} e^{C|s|^\delta} (1+|s|)^{\rho|\si|+B} \sum_{k\leq |s|/\lambda} \frac{(C|s|^{\rho\lambda})^k}{k!} \leq A^{|\si|} (1+|s|)^{\rho|\si|+B} e^{C(|s|^\delta+|s|^{\lambda \rho})} \\
&\ll A^{|\si|} (1+|s|)^{\rho|\si|+B} e^{C|s|^\delta}
\end{split}
\]
with a suitable $\delta <1$, and
\[
S_2 \ll A^{|\si|} \sum_{k> |s|/\lambda} \frac{e^k}{k^k} A^k k^{\rho\lambda k} e^{Ck^\delta} \ll A^{|\si|} \sum_{k> |s|/\lambda}  \frac{A^k}{k^{(1-\lambda \rho)k}} \ll A^{|\si|}
\]
and the result follows.
 
\newpage

\ifx\undefined\bysame{poly}.
\newcommand{\bysame}{\leavevmode\hbox to3em{\hrulefill}\ ,}
\fi

\vskip 1cm
\noindent
Jerzy Kaczorowski, Faculty of Mathematics and Computer Science, A.Mickiewicz University, 61-614 Pozna\'n, Poland and Institute of Mathematics of the Polish Academy of Sciences, 
00-956 Warsaw, Poland. e-mail: kjerzy@amu.edu.pl

\medskip
\noindent
Alberto Perelli, Dipartimento di Matematica, Universit\`a di Genova, via Dodecaneso 35, 16146 Genova, Italy. e-mail: perelli@dima.unige.it

\end{document}